\documentclass[draft]{siamltex}

\usepackage[T1]{fontenc}
\usepackage[latin1]{inputenc}
\usepackage{amsmath}
\usepackage{amssymb}
\usepackage{stmaryrd}
\usepackage{amsfonts}

\newcommand{\q} { \boldsymbol{ q } }
\newcommand{\s} { \boldsymbol{ s } }
\newcommand{\mm} { \boldsymbol{ M } }
\newcommand{\nn} { \boldsymbol{ n } }
\newcommand{\rr} { \boldsymbol{ r } }

\newcommand{\scalpro} { \cdot }
\newcommand{\he} { h_E }
\newcommand{\edge}{E }

\newcommand{\KE}{K(E)}
\newcommand{\mns} {  m _{ns} }

\newcommand{\NABLA} { \boldsymbol{\nabla} }

\newcommand{\eeq} {\end{eqnarray}}

\newcommand{\shear} { \boldsymbol{ Q} }

\newcommand{\vertexset} { \cal{V} }
\newcommand{\ver} { x}
\newcommand{\II} { \boldsymbol{ I } }
\newcommand{\LL} { \boldsymbol{ L } }
\newcommand{\mnn} { m_{nn} }

\newcommand{\MM} { \boldsymbol{ m } }
\newcommand{\MMb} { \boldsymbol{ m }(\boldsymbol{ \beta }) }

\newcommand{\VV} { \boldsymbol{ V } }
\newcommand{\Vstar} { \boldsymbol{ V }\!_* }
\newcommand{\VVh} { \boldsymbol{ V }\!_h }

\newcommand{\R} { \mathbb{ R } }

\newcommand{\nol} { \boldsymbol{ 0 } }

\newcommand{\trinorm}{\Vert\vert}

\newcommand{\bbb} { \boldsymbol{ \beta } }
\newcommand{\yyy} { \boldsymbol{ \eta } }
\newcommand{\fff} { \boldsymbol{ \phi } }
\newcommand{\eee} { \boldsymbol{ \varepsilon } }

\newcommand{\ddiv} { \mathrm{div \,} }
\newcommand{\DDiv} { \mathbf{div \,} }
\newcommand{\la} { \langle }
\newcommand{\ra} { \rangle }
\newcommand{\mc}[1] { \mathcal{ #1 } }

\newcommand{\fa} { \; \; \forall }
\newcommand{\st} { \; \; | \; \; }

\newcommand{\FREE} {\Gamma_{\rm F}}
\newcommand{\FREEh} { \mathcal{F}_h}
\newcommand{\SIMP} {\Gamma_{\rm S}}
\newcommand{\CLAM} {\Gamma_{\rm C}}

\newcommand{\halfFREEh} { \mathcal{F}_{h/2}}

\newcommand{\intedges} { \mathcal{I}_{h}}

\newcommand{\simpedges} { \mathcal{S}_{h}}

\newtheorem{meth}{Method}[section]
\newtheorem{remark}{Remark}[section]
\newtheorem{assu}{Assumption}[section]

\title{A family of ${C}^0$ finite elements for \\
Kirchhoff plates I: Error analysis}

\author{L. Beir${\rm \tilde{a}}$o da Veiga\thanks{
Dipartimento di Matematica "F. Enriques", Universit\`a di Milano, 
via Saldini 50, 20133 Milano, Italy ({\tt beirao@mat.unimi.it})} 
\and J. Niiranen\thanks{
Institute of Mathematics, Helsinki University of Technology, 
P. O. Box 1100, 02015 TKK, Finland ({\tt jarkko.niiranen@tkk.fi})} 
\and R. Stenberg\thanks{
Institute of Mathematics, Helsinki University of Technology, 
P. O. Box 1100, 02015 TKK, Finland ({\tt rolf.stenberg@tkk.fi})}}

\begin{document}

\maketitle


\begin{abstract}
A new finite element formulation for the Kirchhoff plate model is
presented. The method is a displacement formulation with the deflection
and the rotation vector as unknowns and it is based on ideas stemming from
a stabilized method for the Reissner--Mindlin model \cite{S} and a
method to treat a free boundary \cite{DN88}. Optimal a-priori and
a-posteriori error estimates are derived.
\end{abstract}



\pagestyle{myheadings}
\thispagestyle{plain}
\markboth{L. Beir${\rm \tilde{a}}$o da Veiga, J. Niiranen and R. Stenberg}
{A family of ${C}^0$ finite elements for Kirchhoff plates I: Error analysis}

\section{Introduction}
\setcounter{equation}{0}

A conforming finite element method for the Kirchhoff plate bending
problem requires  a $C^1$-continuity and hence leads to methods that
are rarely used in practice. Instead, either a nonconforming method
is used or then the model is abandoned in favor of the
Reissner--Mindlin model. For the latter, there exist several families
of methods that have rigorously been shown to be free from locking
and optimally convergent.

A natural idea is to consider the Kirchhoff model as the limit of
the Reissner--Mindlin model when the plate thickness approaches zero
and to use a good Reissner--Mindlin element with the thickness
(after a scaling, see below) representing the parameter
penalizing the Kirchhoff constraint. In this approach, there are two
obstacles. First, for a free boundary, this leads to a method which
is not consistent. In the literature, this point is often ignored
since mostly the clamped case is considered.
A remedy to this was developed by Destuynder and Nevers who showed,
that the consistency is obtained by adding a term penalizing the
tangential Kirchhoff condition along the free boundary \cite{DN88}.
Even if this modification has been done, there remains a second
drawback. In order that the solution to the penalized formulation is
close to the exact solution, the penalty parameter should be large.
This, however, leads to an ill-conditioned discrete system.

Our aim in the present paper is to present a family of Kirchhoff
plate bending elements  for which   the convergence rate is optimal
even in the presence of free boundaries. The method is a formulation
combining the ideas from the stabilized method for the
Reissner--Mindlin plates presented in \cite{S} and the treatment of
the free boundary presented in \cite{DN88} (cf. \cite{Be,BB} as
well). The family includes "simple low-order" elements and it is
well-conditioned. In the second part \cite{BNS'} of this paper, we
give the results of numerical tests.

The paper is organized as follows. In the next section, we describe
the plate bending problem, and in Section \ref{2}, we introduce the
new family of finite elements. In Section \ref{3}, an a-priori error
analysis is derived. This analysis leads to optimal results, both
with respect to the regularity of the solution and to the polynomial
degree used. In Section \ref{4}, an a-posteriori error analysis is
performed. We derive a local error indicator which is shown to be
both reliable and efficient.

\section{The Kirchhoff plate bending problem} 
\label{1}
\setcounter{equation}{0}

We consider the problem of bending  of an isotropic linearly elastic
plate and assume that the undeformed plate midsurface is described
by a given convex polygonal domain $\Omega \subset \R^2$. The plate
is considered to be clamped on the part $\CLAM$ of its boundary
$\partial\Omega$, simply supported on the part $\SIMP \subset
\partial\Omega$ and free on $\FREE \subset
\partial\Omega$. The deflection and transversal load are  denoted by $w$ and $g$,
respectively.

In the sequel, we indicate with $\vertexset$ the set of all corner
points in $\FREE$. Moreover, $\nn$ and $\s$ represent the unit
outward normal and the unit counterclockwise tangent to the
boundary. Finally, for points $\ver \in \vertexset$, we introduce
the following notation. We indicate with $\nn_1$ and $\s_1$ the unit
vectors corresponding, respectively, to $\nn$ and $\s$ on one of the
two edges   forming the boundary angle at $\ver$;  with $\nn_2$ and
$\s_2$ we indicate the ones corresponding to the other edge. Note
that which of the two edges correspond to the subscript $1$ or $2$
is not relevant.

The classical Kirchhoff plate bending model is then given by the
biharmonic partial differential equation
\begin{equation}
\label{B-pbl}
  \mathsf{D} \Delta^2 w =  g
\qquad  {\rm in} \ \Omega \, ,
\end{equation}
the  boundary conditions
\begin{equation}
\begin{array}{llll}
   w =0 \, , &    \ \ \frac{\partial w}{\partial \nn}  = 0
& & {\rm on} \ \CLAM \, ,
  \\
& &
 \\
  w =0 \, , &  \ \ \nn\scalpro \mm \nn  = 0
& & {\rm on} \ \SIMP \, ,
   \\
& & \\
\nn\scalpro \mm \nn = 0 \, ,&    \ \ \frac{\partial}{\partial \s}
\big(\s\scalpro \mm \nn \big)+ (\DDiv \mm)\cdot\nn  = 0
& & {\rm on} \ \FREE \, ,
\end{array}
\end{equation}
and the corner conditions
\begin{equation}
\label{cornercond}
\big(\s_1\scalpro \mm\nn_1 \big)(\ver) = \big(\s_2\scalpro \mm
\nn_2\big) (\ver) \qquad \forall \ver \in \vertexset \, .
\end{equation}
Here
\begin{equation}
\mathsf{D} = \frac{\mathsf{E}t^3}{12(1-\nu^2)}
\end{equation}
is the bending rigidity,  with $ \mathsf{E} $, $\nu$ being the Young modulus
and the Poisson ratio for the material, respectively.  Note that for
the shear modulus  $ \mathsf{G} $ it holds
\begin{equation}
\mathsf{G} = \frac{\mathsf{E}}{2(1+\nu)} \, .
\end{equation}
The moment tensor is given by
\begin{equation}
\label{omomdef}
  \mm(\nabla w) = \mathsf{D} \big( (1-\nu) \eee(\nabla w) 
                  + \nu \, \ddiv\! (\nabla w )\, \II \big) \, ,
\end{equation}
with the symmetric gradient $ \eee $, and the shear force by
\begin{equation}
\label{oshear}
\shear =-\DDiv \mm \, .
\end{equation}
Note, that the independence of the Poisson ratio $\nu$ in the
differential equation (\ref{B-pbl}) is a consequence of cancellations
when substituting (\ref{omomdef}) and (\ref{oshear}) into the
equilibrium equation
\begin{equation}
-\ddiv \shear = g \, .
\end{equation}
For the analysis below, it will be convenient to perform a scaling
of the problem by assuming that the load is given by $g= \mathsf{G} t^3 f$,
with $f$ fixed. Then the differential equation (\ref{B-pbl}) becomes
independent of the plate thickness:
\begin{equation}
\label{scaledbi}
 \frac{1}{ 6(1 - \nu)} \, \Delta^2 w =  f,\ \qquad  {\rm in} \ \Omega \, .
\end{equation}
Furthermore, we use the following scaled moment tensor $\MM$ given
by
\begin{equation}
\mm(\nabla w)= \mathsf{G} t^3 \MM(\nabla w)
\end{equation}
and the shear force $\q$ is defined by
\begin{equation}
\shear = \mathsf{G} t^3 \q \, .
\end{equation}

The unknowns in our finite element method will be the approximations to
the deflection and its gradient, the rotation $\bbb=\nabla w$. With
this as a new unknown, our problem can be written as the system of
partial differential equations
\begin{eqnarray}
 \nabla w - \bbb &= &\nol \, ,
\label{strong3}
 \\
\label{strong2}
     -  \ddiv \q  &=& f \, ,
  \\
\LL \bbb + \q  &=& \nol \, , 
\qquad \textrm{in } \Omega \, ,
\label{strong1}
\end{eqnarray}
the boundary conditions
\begin{equation}
\label{bc-clamed}
 w=0 \, , \ \bbb = \nol \, ,  
 \qquad {\rm on} \ \CLAM \, ,  \\
\end{equation}
\begin{equation}
\label{bc-simp}
 w=0 \, , \ \bbb\cdot\s = 0 \, , \ \nn\scalpro\MMb \nn = 0 \, ,
\qquad {\rm on} \ \SIMP \, ,
\end{equation}
\begin{equation}
\label{bc-free}
\frac{\partial w}{\partial \s} - \bbb\cdot\s=0 \, , \
\nn\scalpro\MMb \nn = 0  \, , \ \frac{\partial}{\partial \s}
\big(\s\scalpro\MMb
\nn\big) - \q\cdot\nn = 0 \, , \  {\rm on} \ \FREE \, , \\
\end{equation}
and the corner conditions
\begin{equation}
\label{corner}
 \big(\s_1\scalpro\MMb \nn_1 \big)(\ver) = \big(\s_2\scalpro\MMb \nn_2 \big)(\ver) 
 \qquad \forall \ver \in \vertexset \, .
\end{equation}
 The operator $ \LL $ is defined as
\begin{equation}
\label{L-op}
\LL \bbb= \DDiv\MMb   \, ,
\end{equation}
and the  scaled bending moment is considered as a function of the
rotation:
\begin{equation}
\label{BM}
\MM(\bbb)=\frac{1}{6} \big( \eee(\bbb) 
+ \frac{\nu}{1 - \nu} \ddiv \bbb \, \II \big) \, .
\end{equation}
In the sequel,  we will often write $\MM$ instead of $\MM(\bbb)$. We
further denote
\begin{equation}
 \label{adef}
a(\bbb,\yyy)= (\MM(\bbb),\eee(\yyy)) \, .
\end{equation}

In order to neglect plate rigid movements and the related
technicalities, we will in the sequel assume that the
one-dimensional measure of $ \CLAM  $ is positive.

\section{The finite element formulation}
\label{2}
\setcounter{equation}{0}

In this section, we will
introduce our finite element method. Even if our method  is stable
for all choices of finite element spaces, we will, for simplicity,
present it for triangular elements and for the polynomial degrees
that yield an optimal convergence rate. Hence,
 let a regular family of triangular meshes on $\Omega$ be given.
For the  integer $k\ge 1$, we then define the discrete spaces
\begin{align}
\label{disc-spaces}
    W_h & = \{ w \in W \st w_{|K} \in P_{k+1}(K) \fa K \in \mc{ C }_h \} \, ,
\\
    \VVh & = \{ \yyy \in \VV \st \yyy_{|K} \in [P_k(K)]^2 \fa K \in \mc{ C }_h \} \, ,
\end{align}
with
\begin{eqnarray}
\label{Wspace}
W &=& \left\{ v \in H^1(\Omega) \: \vert \: 
v=0 \ {\rm on} \ \CLAM \cup \SIMP \right\} \, , 
\\
\label{Vspace} 
\VV &=& \left\{ \yyy \in [H^1(\Omega)]^2 \: \vert \:
\yyy=\nol \ {\rm on} \ \CLAM \: , \ \yyy\cdot\s=0 \ {\rm on} \ \SIMP \right\} \, .
\end{eqnarray}
Here $\mc{C}_h$ represents the set of all triangles $K$ of the
mesh and $P_k(K)$ is the space of polynomials of degree $k$ on $K$.
In the sequel, we will indicate with $h_K$ the diameter of each
element $K$, while $h$ will indicate the maximum size of all the
elements in the mesh. Furthermore, we will indicate with $\edge$  a
general edge of the triangulation and with $\he$ the length of
$\edge$. The set of all edges lying on the free boundary $\FREE$ we
denote by $\FREEh$.

Before introducing the method, we state the following result which
trivially follows  from classical scaling arguments and the
coercivity of the form $a$.

\begin{lemma}
\label{inverse}
There exist positive constants $C_I$ and $C_I'$ such that
\begin{eqnarray}
&& C_I \sum_{K \in \mc{C}_h} h_K^2 \Vert \LL \fff \Vert_{0,K}^2 \le a(\fff,\fff) 
\qquad \forall \fff \in \VVh \, , 
\\
&& C_I'  \sum_{\edge \in \FREEh} \!\!\! \he \: \Vert \mns(\fff)
\Vert_{0, \edge }^2 \le a(\fff,\fff) \qquad \forall \fff \in \VVh \, ,
\end{eqnarray}
where the operator $\mns(\fff) = \s\scalpro \MM (\fff)   \nn$ with
$\nn,\s$, being the unit outward normal and the unit
counterclockwise tangent to the edge $ \edge $, and with $\MM$
defined in $(\ref{BM})$.
\end{lemma}

Let two real numbers $ \gamma $ and $ \alpha $ be assigned, $ \gamma
> 2/C_I' $ and $ 0 < \alpha < C_I/4 $. Then, the discrete problem
reads as follows.

\begin{meth}
\label{meth} 
Find $ (w_h, \bbb_h) \in W_h \times \VVh $, such that
\begin{equation}
\label{metheq}
    \mc{A}_h(w_h, \bbb_h; v, \yyy) = (f, v)  \quad \fa (v, \yyy) \in W_h \times \VVh \, ,
\end{equation}
where the form $ \mc{A}_h $ is defined as
\begin{eqnarray}
\mc{A}_h (z, \fff; v, \yyy) = \mc{B}_h (z, \fff; v, \yyy) 
+ \mc{D}_h (z,\fff; v, \yyy) \, ,
\end{eqnarray}
with
\begin{eqnarray}
\label{added2}
 \mc{B}_h (z, \fff; v, \yyy) 
& = & a(\fff, \yyy) - \sum_{K \in \mc{ C }_h} \alpha h_K^2 (\LL \fff, \LL \yyy)_K
\nonumber \\
&& + \sum_{K \in \mc{ C }_h} \frac{1}{\alpha h_K^2} (\nabla z
- \fff - \alpha h_K^2 \LL \fff, \nabla v - \yyy - \alpha h_K^2 \LL \yyy)_K 
\end{eqnarray}
and
\begin{eqnarray}
\label{Ddef}
  \mc{D}_h (z,\fff; v, \yyy)& =  &   \langle \mns(\fff),[\nabla v -
\yyy]\cdot \s\rangle_{\FREE} + \langle [\nabla z - \fff]\cdot\s,
\mns(\yyy)\rangle_{\FREE}
\nonumber \\
&&  + \sum_{\edge \in \FREEh }\frac{\gamma}{\he} \langle[\nabla z -
\fff]\cdot\s , [\nabla v - \yyy]\cdot \s\rangle_\edge
\end{eqnarray} 
for all $(z, \fff), (v, \yyy) \in W_h \times \VVh$.
Here $\langle\cdot, \cdot\rangle_{\FREE}$ and $\langle \cdot,
\cdot\rangle_E$ denote the $L^2$-inner products on $\FREE$ and $E$,
respectively.
\end{meth}

The bilinear form $\mc{B}_h$ constitutes the Reissner--Mindlin method
of \cite{S} with the thickness $t$ set equal to zero, while the
additional form $\mc{D}_h$ is introduced in order to avoid the convergence
deterioration in the presence of free boundaries.

Furthermore, we introduce the discrete shear force
\begin{equation}
\label{disc-shear}
\q_{h|K} = \frac{1}{\alpha h_K^2} (\nabla w_h - \bbb_h - \alpha
h_K^2 \LL \bbb_h)_{ |K } \fa K \in \mc{ C }_h \, .
\end{equation}
We note that, due to (\ref{strong1}) and (\ref{strong3}), it holds
\begin{equation}
\label{shear}
\q_{\vert K} = \frac{1}{\alpha h_K^2} (\nabla w - \bbb - \alpha
h_K^2 \LL \bbb)_{ |K } \fa K \in \mc{ C }_h \, ,
\end{equation}
and hence it follows that the definition (\ref{disc-shear}) is
consistent with the exact shear force.

For simplicity, in the rest of this section we assume that the
deflection $w$   belongs to  $H^3(\Omega)$; this is a very
reasonable assumption, as discussed at the end of this section. Note
as well that, with some additional technical work involving the
appropriate Sobolev spaces and their duals, such assumption could be
probably avoided. The following result states the consistency of the
method.

\begin{theorem}
\label{consistency}
The solution $(w,\bbb)$ of the problem
$(\ref{strong1})$--$(\ref{corner})$ satisfies
\begin{equation}
\mc{A}_h(w, \bbb; v, \yyy) = (f, v)  
\quad \fa (v, \yyy) \in W_h \times \VVh \, .
 \end{equation}
\end{theorem}
\begin{proof}
The definition of the bilinear forms in Method \ref{meth}, recalling
(\ref{strong1}) and the expression (\ref{shear}), give
\begin{eqnarray}
\label{first}
&& \mc{B}_h (w, \bbb; v, \yyy) = a(\bbb, \yyy) - \sum_{K \in \mc{ C
}_h} \alpha h_K^2 (\LL \bbb, \LL \yyy)_K
\nonumber \\
&& \qquad + \sum_{K \in \mc{ C }_h} \frac{1}{\alpha h_K^2} (\nabla w
- \bbb - \alpha h_K^2 \LL \bbb, \nabla v - \yyy - \alpha h_K^2 \LL
\yyy)_K
\nonumber \\
&& \quad = a(\bbb, \yyy) + \sum_{K \in \mc{ C }_h} \alpha h_K^2 (\q, \LL
\yyy)_K + \sum_{K \in \mc{ C }_h} (\q, \nabla v - \yyy - \alpha
h_K^2 \LL \yyy)_K
\nonumber \\
&& \quad = a(\bbb, \yyy) + (\q, \nabla v - \yyy) \, .
\end{eqnarray}
First, by the definition (\ref{adef}), then integrating by parts on
each triangle, finally using the regularity of the functions
involved, and the boundary conditions (\ref{bc-clamed}),
(\ref{bc-simp})  on $\CLAM, \ \SIMP$, respectively, we get
\begin{eqnarray}
&& a(\bbb,\yyy) + (\q,\nabla v -\yyy) =
(\MM(\bbb),\eee(\yyy))  + (\q,\nabla v -\yyy)
\nonumber \\
& & \quad = - (\LL\bbb + \q, \yyy) + \langle
\MM(\bbb)\cdot\nn,\yyy\rangle_{\FREE} - (\ddiv \q,v) + \langle
\q\cdot\nn,v\rangle_{\FREE} \, .
\end{eqnarray}
Recalling (\ref{strong1}) and (\ref{strong2}), the identity above
becomes
\begin{eqnarray}
a(\bbb,\yyy) + (\q,\nabla v -\yyy) = (f,v) + \langle
\MM(\bbb)\cdot\nn,\yyy\rangle_{\FREE} + \langle
\q\cdot\nn,v\rangle_{\FREE} \, ,
\end{eqnarray}
while, using the boundary conditions of (\ref{bc-free}) on $\FREE$
and integration by parts along the boundary, finally leads to
\begin{eqnarray}
\label{first-bis}
a(\bbb,\yyy) + (\q,\nabla v -\yyy) = (f,v) - \langle
m_{ns}(\bbb),[\nabla v - \yyy]\cdot\s\rangle_{\FREE} \, .
\end{eqnarray}
Due to (\ref{bc-free}), we have
\begin{eqnarray}
\label{second}
 \mc{D}_h (w,\bbb; v, \yyy) &= & \langle m_{ns}(\bbb),[\nabla v -
\yyy]\cdot \s\rangle_{\FREE} + \langle[\nabla w - \bbb]\cdot\s,
m_{ns}(\yyy)\rangle_{\FREE}
\nonumber \\
&& \quad +  \sum_{\edge \in \FREEh}\frac{\gamma}{h_\edge} \langle
[\nabla w - \bbb]\cdot\s , [\nabla v - \yyy]\cdot \s\rangle_\edge
\nonumber \\
&  =&   \langle m_{ns}(\bbb),[\nabla v - \yyy]\cdot
\s\rangle_{\FREE}  \, .
\end{eqnarray}
The result now directly follows from (\ref{first}),
(\ref{first-bis}) and (\ref{second}).
\end{proof}

\begin{remark}
\label{rem1}
If the Reissner--Mindlin method of \cite{S} without the
additional form $ \mc{D}_h $ is employed by setting $t=0$, then in
the presence of a free boundary we obtain
\begin{equation}
\label{incons}
\mc{B}_h(w, \bbb; v, \yyy) = (f, v) +  \langle \mns(\bbb),[\nabla v
- \yyy]\cdot \s\rangle_{\FREE}\qquad \fa (v, \yyy) \in W_h \times \VVh \, .
\end{equation}
Therefore, this would lead to an inconsistent method.  We return to
this in Remark  \ref{rem2} below.
\end{remark}

\section{Stability and a-priori error estimates}
\label{3}
\setcounter{equation}{0}

For $(v,\yyy) \in W_h \times \VVh$, we introduce the following
mesh dependent norms:
\begin{eqnarray}
\label{eq:h-norm}
&& \vert (v,\yyy) \vert_h^2 = \sum_{K\in\mc{C}_h} h_K^{-2} 
\Vert \nabla v -\yyy \Vert_{0,K}^2 \, ,  \\
\label{2hnorm} 
&& \Vert v \Vert_{2,h}^2 = \Vert v \Vert_{1}^2 +
\sum_{K\in\mc{C}_h} \vert v \vert_{2,K}^2 + \sum_{\edge\in\intedges}
h_\edge^{-1} \Vert \, \llbracket \frac{\partial v}{\partial \nn}
\rrbracket \, \Vert_{0,\edge}^2 +  \sum_{\edge\subset \CLAM}
h_\edge^{-1} \Vert \frac{\partial v}{\partial \nn} \Vert_{0,\edge}^2 \, , 
\\
\label{3bar-norm} 
&& \trinorm (v,\yyy) \trinorm_h = \Vert \yyy
\Vert_1 + \Vert v \Vert_{2,h} + \vert (v,\yyy) \vert_h \, ,
\end{eqnarray}
where $ \llbracket \cdot \rrbracket$ represents the jump operator
and $ \intedges $ denotes the edges lying in the interior 
of the domain $ \Omega $.

In \cite{Pi88}, the following lemma is proved.

\begin{lemma}
There exists a positive constant $C$ such that 
\begin{eqnarray} 
\Vert v \Vert_{2,h} 
\le C \big( \Vert \yyy \Vert_1 + \Vert v \Vert_1  
+ \vert (v,\yyy) \vert_h \big) 
\qquad \forall (v,\yyy) \in W_h \times \VVh  \, . 
\end{eqnarray}
\end{lemma}

Using the Poincar\'e inequality and the previous lemma, the
following equivalence easily follows:

\begin{lemma}
\label{equivalence}
There exists a positive constant $C$ such that
\begin{equation}
C \trinorm (v,\yyy) \trinorm_h \le \Vert \yyy \Vert_1 + \vert
(v,\yyy) \vert_h \le \trinorm (v,\yyy) \trinorm_h \quad \forall
(v,\yyy) \in W_h \times \VVh \, .
\end{equation}
\end{lemma}

We now have the following stability estimate.

\begin{theorem}
\label{stability}
Let $ 0 < \alpha < C_I/4 $ and $ \gamma > 2/C_I' $.
Then there exists a positive constant $C$ such that
\begin{equation}
\mc{A}_h(v, \yyy; v, \yyy) \ge C \trinorm (v,\yyy) \trinorm_h^2
\quad \forall (v,\yyy) \in W_h \times \VVh \, .
\end{equation}
\end{theorem}

\begin{proof}
Using the first inverse estimate of Lemma \ref{inverse} we  get
\begin{eqnarray}
\label{eq1}
&& \mc{B}_h(v, \yyy; v, \yyy)
\nonumber \\
&& \quad = a(\yyy,\yyy) - \sum_{K \in \mc{ C }_h} \alpha h_K^2 \Vert \LL \yyy \Vert_{0,K}^2
+ \sum_{K \in \mc{ C }_h} \frac{1}{\alpha h_K^2} 
\Vert \nabla v - \yyy - \alpha h_K^2 \LL \yyy \Vert_{0,K}^2 
\nonumber \\
&& \quad \ge \big( 1-\frac{\alpha}{C_I} \big) \: a(\yyy,\yyy) 
+ \sum_{K \in \mc{ C }_h} \frac{1}{\alpha h_K^2} \Vert \nabla v - \yyy - \alpha h_K^2 \LL \yyy \Vert_{0,K}^2 \, .
\end{eqnarray}
Next,  using locally the arithmetic-geometric mean inequality with
the  constant $\gamma/\he$, then the second inverse inequality of
Lemma \ref{inverse}, we get
\begin{eqnarray}
\label{Dh}
&& \mc{D}_h(v, \yyy; v, \yyy)
\nonumber \\
&& \quad = \sum_{\edge \in \FREEh} \Big( 2 \langle \mns(\yyy),[\nabla v -
\yyy]\cdot \s\rangle_\edge + \frac{\gamma}{\he} \Vert [\nabla v -
\yyy]\cdot\s \Vert_{0,\edge}^2 \Big)
\nonumber \\
&& \quad \ge \sum_{\edge \in \FREEh} \Big( - \frac{\gamma}{\he} \Vert
[\nabla v - \yyy]\cdot\s \Vert_{0,\edge}^2 -  \gamma^{-1} \he  \:
\Vert \mns(\yyy) \Vert_{0,\edge}^2 + \frac{\gamma}{\he} \Vert
[\nabla v - \yyy]\cdot\s \Vert_{0,\edge}^2 \Big)
\nonumber \\
&& \quad  = - \sum_{\edge \in \FREEh} \gamma^{-1} \he  \: \Vert
\mns(\yyy) \Vert_{0,\edge}^2
\nonumber \\
&& \quad \ge - \frac{\gamma^{-1}}{C_I'} \: a(\yyy,\yyy)
\ge - \frac{1}{2} \: a(\yyy,\yyy) \, .
\end{eqnarray}
Joining (\ref{eq1}) with (\ref{Dh}) and using Korn's inequality we then obtain
\begin{eqnarray}
\label{sum}
&& \mc{B}_h(v, \yyy; v, \yyy) + \mc{D}_h(v, \yyy; v, \yyy)
\nonumber \\
&& \quad \ge \big( \frac{1}{2}-\frac{\alpha}{C_I} \big) \: a(\yyy,\yyy)
+ \sum_{K \in \mc{ C }_h} \frac{1}{\alpha h_K^2} \Vert \nabla v - \yyy - \alpha h_K^2 \LL \yyy \Vert_{0,K}^2
\nonumber \\
&& \quad \ge C \Big( \Vert \yyy \Vert_1^2
+ \sum_{K \in \mc{ C }_h} \frac{1}{\alpha h_K^2} \Vert \nabla v - \yyy - \alpha h_K^2 \LL \yyy \Vert_{0,K}^2 \Big) \, .
\end{eqnarray}
From the triangle inequality, again the inverse estimate of Lemma
\ref{inverse} and the boundedness of the bilinear form $a$ it
follows
\begin{eqnarray}
&& \sum_{K \in \mc{ C }_h} \frac{1}{\alpha h_K^2} \Vert \nabla v - \yyy \Vert_{0,K}^2
\nonumber \\
&& \quad \le 2 \Big( \sum_{K \in \mc{ C }_h} \frac{1}{\alpha h_K^2} 
\Vert \nabla v - \yyy - \alpha h_K^2 \LL \yyy \Vert_{0,K}^2 
+ \sum_{K \in \mc{ C }_h} \frac{1}{\alpha h_K^2} \Vert \alpha h_K^2 \LL \yyy \Vert_{0,K}^2 \Big)
\nonumber \\
&& \quad \le 2 \Big( \sum_{K \in \mc{ C }_h} \frac{1}{\alpha h_K^2} 
\Vert \nabla v - \yyy - \alpha h_K^2 \LL \yyy \Vert_{0,K}^2 
+ \sum_{K \in \mc{ C }_h} \alpha h_K^2 \Vert \LL \yyy \Vert_{0,K}^2 \Big)
\nonumber \\
&& \quad \le C \Big( \sum_{K \in \mc{ C }_h} \frac{1}{\alpha h_K^2} 
\Vert \nabla v - \yyy - \alpha h_K^2 \LL \yyy \Vert_{0,K}^2 + a(\yyy,\yyy) \Big)
\nonumber \\
&& \quad \le C \Big( \sum_{K \in \mc{ C }_h} \frac{1}{\alpha h_K^2} 
\Vert \nabla v - \yyy - \alpha h_K^2 \LL \yyy \Vert_{0,K}^2 
+ \Vert \yyy \Vert_1^2 \Big) \, ,
\end{eqnarray}
which combined with (\ref{sum}) gives
\begin{equation}
\mc{A}_h(v, \yyy; v, \yyy) \ge C \big( \Vert \yyy \Vert_1^2 + \vert
(v,\yyy) \vert_h ^2 \big) \, .
\end{equation}
The result then follows from the norm equivalence of Lemma
\ref{equivalence}.
\end{proof}

We can now derive the error estimates for the method.
We note that the assumptions of the theorem are
supposed to be valid for the further results below as well,
hence not repeated in the sequel.

\begin{theorem}
\label{apriori}
Let $0 < \alpha < C_I/4 $ and $ \gamma > 2/C_I' $. Let $(w,\bbb)$ be
the exact solution of the problem and let $(w_h,\bbb_h)$ be the
approximate solution obtained with Method \ref{meth}. Suppose that
$w\in H^{s+2}(\Omega)$, with $1\leq s \leq k$. Then it holds
\begin{eqnarray}
\label{added}
\trinorm (w-w_h,\bbb-\bbb_h) \trinorm_h \le C h^s \Vert w \Vert_{s+2} \, .
\end{eqnarray}
\end{theorem}

\begin{proof}
{\em Step 1.} Let $(w_I,\bbb_I) \in W_h \times \VVh$ be the usual
Lagrange interpolants to $w$ and $\bbb$, respectively. Using first
the stability result of Theorem \ref{stability} and then the
consistency result of Theorem \ref{consistency} one has the
existence of a pair
\begin{equation}
\label{1-bound}
(v,\yyy) \in W_h\times\VVh \, , 
\qquad \trinorm (v,\yyy) \trinorm_h \le C \, ,
\end{equation}
such that
\begin{eqnarray}
\label{start}
\trinorm (w_h-w_I,\bbb_h-\bbb_I) \trinorm_h
& \le & \mc{A}_h(w_h-w_I, \bbb_h-\bbb_I; v, \yyy)
\nonumber \\
& = & \mc{A}_h(w-w_I, \bbb-\bbb_I; v, \yyy) \, ,
\end{eqnarray}
where we recall that $ \mc{A}_h = \mc{B}_h + \mc{D}_h $.

{ \em Step 2.} For the $ \mc{B}_h $-part, we have
\begin{eqnarray}
\label{Bhbound}
&& \mc{B}_h(w-w_I, \bbb-\bbb_I; v, \yyy)
 = a(\bbb-\bbb_I, \yyy)
 - \sum_{K \in \mc{ C }_h} \alpha h_K^2 (\LL (\bbb - \bbb_I), \LL \yyy)_K
\\ \nonumber
&& \quad + \sum_{K \in \mc{ C }_h} \frac{1}{\alpha h_K^2} (\nabla(w-w_I) - (\bbb-\bbb_I)
   - \alpha h_K^2 \LL(\bbb -\bbb_I),
 \nabla v - \yyy - \alpha h_K^2 \LL \yyy)_K \, .
\end{eqnarray}
Due to the first inverse inequality of Lemma \ref{inverse}, we get
\begin{equation}
\Big( \sum_{K \in \mc{C}_h} h_K^2 \Vert \LL \yyy \Vert_{0,K}^2 \Big)^{1/2} 
\le C \trinorm (v,\yyy) \trinorm_h
\end{equation}
and
\begin{equation}
\Big( \sum_{K \in \mc{C}_h} \frac{1}{\alpha h_K^2} 
\Vert \nabla v - \yyy - \alpha h_K^2 \LL \yyy \Vert_{0,K}^2 \Big)^{1/2}
\le C \trinorm (v,\yyy) \trinorm_h \, .
\end{equation}
Using these bounds in (\ref{Bhbound}) and recalling (\ref{1-bound}) we obtain
\begin{eqnarray}
\label{x1}
&& \mc{B}_h(w-w_I, \bbb-\bbb_I; v, \yyy)
\nonumber \\
&& \quad \le C \Big (\trinorm (w-w_I,\bbb-\bbb_I)\trinorm_h 
+ \big( \sum_{K \in \mc{C}_h} h_K^2 \vert \bbb - \bbb_I \vert_{2,K}^2 \big)^{1/2} \Big) \, .
\end{eqnarray}
Substituting the definition of the norm (\ref{3bar-norm}) in (\ref{x1}), using the triangle inequality, and finally applying the classical interpolation estimates it easily follows
\begin{eqnarray}
\label{B-end}
\mc{B}_h(w-w_I, \bbb-\bbb_I; v, \yyy) 
\le C h^s \big( \Vert w \Vert_{s+2} + \Vert \bbb \Vert_{s+1} \big) \, .
\end{eqnarray}

{\em Step 3.} For the $ \mc{D}_h $-part in (\ref{start}), we have,
by the definition (\ref{Ddef}),
\begin{eqnarray}
\label{D-start}
&& \mc{D}_h (w-w_I,\bbb-\bbb_I; v, \yyy) =
\langle\mns(\bbb-\bbb_I),[\nabla v - \yyy]\cdot \s\rangle_{\FREE}
\nonumber \\
&& \qquad + \langle[\nabla (w-w_I) - (\bbb-\bbb_I)]\cdot\s,
\mns(\yyy)\rangle_{\FREE}
\nonumber \\
&& \qquad + \sum_{\edge \in \FREEh}\frac{\gamma}{\he} \langle[\nabla
(w-w_I) - (\bbb-\bbb_I)]\cdot\s , [\nabla v - \yyy]\cdot
\s\rangle_\edge
\nonumber \\
&& \quad =: T_1 + T_2 + T_3 \, .
\end{eqnarray}
Scaling arguments give
\begin{eqnarray}
\label{x2} 
\Vert [\nabla v - \yyy]\cdot \s \Vert_{0,\edge}^2 \le
\Vert \nabla v - \yyy \Vert_{0,\edge}^2 \le C h_{\KE}^{-1} \Vert
\nabla v - \yyy \Vert_{0,\KE}^2
\end{eqnarray}
for all $\edge \in \FREEh$, where $\KE$ is the
triangle with $\edge$ as an edge. The $l^2$-Cauchy--Schwartz inequality,
the bound (\ref{x2}) and  the norm definition (\ref{3bar-norm}) now
give
\begin{eqnarray} 
T_1 & \le & \Big( \sum_{\edge \in \FREEh}
h_{\KE} \Vert \mns(\bbb-\bbb_I) \Vert_{0,\edge}^2 \Big)^{1/2} \Big(
\sum_{\edge \in \FREEh} h_{\KE}^{-1} \Vert [\nabla v - \yyy]\cdot \s
\Vert_{0,\edge}^2 \Big)^{1/2}
\nonumber \\
& \le & C \Big( \sum_{\edge \in \FREEh} h_{\KE} \Vert
\mns(\bbb-\bbb_I) \Vert_{0,\edge}^2 \Big)^{1/2} \trinorm (v,\yyy)
\trinorm_h \, .
\end{eqnarray}
Recalling the bound (\ref{1-bound}), classical
polynomial interpolation properties give
\begin{eqnarray}
\label{T1}
 T_1   \le   C \Big( \sum_{\edge \in \FREEh}
h_{\KE} \Vert \mns(\bbb-\bbb_I) \Vert_{0,\edge}^2 \Big)^{1/2}
 \le  C h^{s} \Vert \bbb \Vert_{s+1} \, .
\end{eqnarray}

Again, by scaling we have
\begin{eqnarray}
\label{x3}
\Vert \mns(\yyy) \Vert_{0,\edge}^2 \le h_{\KE}^{-1} \vert
\yyy \vert_{1,\KE}^2 \qquad \forall \edge \in \FREEh \, .
\end{eqnarray}
The $l^2$-Cauchy--Schwartz inequality, this bound
and the norm definition (\ref{3bar-norm}) give
\begin{eqnarray}
T_2 & \le & \Big( \sum_{\edge \in \FREEh}  h_{\KE}^{-1} \Vert \nabla
(w-w_I) - (\bbb-\bbb_I)\Vert_{0,\edge}^2 \Big)^{1/2} \Big(
\sum_{\edge \in \FREEh} h_{\KE} \Vert \mns(\yyy) \Vert_{0,\edge}^2
\Big)^{1/2}
\nonumber \\
& \le & C \Big( \sum_{\edge \in \FREEh} h_{\KE}^{-1} \Vert  \nabla
(w-w_I) - (\bbb-\bbb_I) \Vert_{0,\edge}^2 \Big)^{1/2} \trinorm
(v,\yyy) \trinorm_h \, .
\end{eqnarray}
Recalling the bound (\ref{1-bound}),   classical polynomial
interpolation estimates give
\begin{eqnarray}
\label{T2} T_2 & \le &
C\Big( \sum_{\edge \in \FREEh} h_{\KE}^{-1} \Vert  \nabla
(w-w_I) - (\bbb-\bbb_I) \Vert_{0,\edge}^2 \Big)^{1/2}
\nonumber \\
& \le & C h^{s} \big( \Vert \bbb \Vert_{s+1} + \Vert w \Vert_{s+2} \big) \, .
\end{eqnarray}

The bound for $T_3$ follows combining the same techniques used for $T_1$ and $T_2$; we get
\begin{eqnarray}
\label{T3}
T_3 \le C h^{s} \big( \Vert \bbb \Vert_{s+1} + \Vert w \Vert_{s+2} \big) \, .
\end{eqnarray}
Now, joining all the bounds (\ref{start}), (\ref{B-end}), (\ref{D-start}), (\ref{T1}), (\ref{T2}) and (\ref{T3}) we obtain
\begin{eqnarray}
\label{almost}
\trinorm (w_h-w_I,\bbb_h-\bbb_I) \trinorm_h \le C h^{s} \big( \Vert \bbb \Vert_{s+1} + \Vert w \Vert_{s+2} \big) \, .
\end{eqnarray}
The triangle inequality and the classical polynomial interpolation estimates (recalling that $\bbb = \nabla w$) then yield
\begin{equation}
\label{final-apriori}
\trinorm (w-w_h,\bbb-\bbb_h) \trinorm_h  
\le C h^{s} \big( \Vert \bbb \Vert_{s+1} + \Vert w \Vert_{s+2} \big)
\le C h^{s} \Vert w \Vert_{s+2} \, .
\end{equation}
Note that the result holds for real values of the regularity parameter $s$ 
since the interpolation results used above are valid for real values of $s$.
\end{proof}

\begin{remark}
\label{rem2} 
As noted in Remark \ref{rem1},
the  limiting Reissner--Mindlin  method (i.e.,
without the additional correction $ \mc{D}_h  $) is inconsistent.
Regardless of the solution regularity and the polynomial degree $k$,
the inconsistency term   can only be bounded with the order
$O(h^{1/2})$. As well known (see for example \cite{PS96}), the
inconsistency error is a lower bound for the error of finite element
methods. As a consequence, the numerical scheme will not converge with a
rate better than $h^{1/2}$ if $\FREE \not= \emptyset$. This
observation is also confirmed by the numerical tests shown in
\cite{BNS'}. See \cite{C} for other numerical tests regarding this
issue. Note further, that this boundary inconsistency term is
connected not only to the formulation in \cite{S} but is common to
any other Kirchhoff method which follows a "Reissner--Mindlin limit"
approach.
\end{remark}

For the shear force, the practical norm to use is the discrete
negative norm
\begin{equation}
\label{shearnorm}
\Vert \rr   \Vert_{-1,h}  
= \big( \sum_{K\in \mc{C}_h} h_K^2\Vert\rr\Vert_{0,K}^2\big)^{1/2} \, .
\end{equation}
Since we assume that $w\in H^{s+2}(\Omega) $, with $s\geq1$, we
have $\q\in [L^2(\Omega)]^2$ and from the estimates above it
immediately follows:

\begin{lemma}
\label{shearlemma2} 
It holds
\begin{eqnarray}
\label{minusoneh}
\Vert \q -\q_h \Vert_{-1,h} \le C h^s \Vert w \Vert_{s+2} \, .
 \end{eqnarray}
\end{lemma}

From this it follows a norm estimate in the dual to the space
\begin{equation}
\label{vstar} 
\Vstar = \left\{ \yyy \in [H^1(\Omega)]^2 \: \vert \:
\yyy= {\bf 0} \ {\rm on} \ \CLAM, \ \yyy\cdot\s=0 \ {\rm on} \ \FREE
\cup \SIMP  \right\} \, ,
\end{equation}
i.e., in the norm
\begin{equation}
\label{eq:neg*norm}
\Vert \rr \Vert_{-1,*} = \sup_{\yyy \in \Vstar}
\frac{\langle\rr,\yyy\rangle}{\Vert \yyy \Vert_1} \, .
\end{equation}
The result we have is the following.

\begin{lemma}
\label{shearlemma3}
It holds 
\begin{equation} 
\Vert \q -\q_h \Vert_{-1,*} \le C h^s \Vert w \Vert_{s+2} \, .
\end{equation}
\end{lemma}

\begin{proof}
The proof is essentially an application of the
"Pitk\"{a}ranta--Verf\"{u}rth trick" (see \cite{Pi79,Ve84}).
By the definition of the norm $\Vert\cdot\Vert_{-1,*}$ there exists a
function $\yyy \in \Vstar$ such that
\begin{equation}
\label{zz-1} 
\Vert \q -\q_h \Vert_{-1,*} \le (\q -\q_h, \yyy) \, ,
\qquad \Vert \yyy \Vert_{1} \le C \, .
\end{equation}
Using a Cl\'ement
type interpolant we can find a piecewise linear function $\yyy_I \in
\Vstar$ such that it holds
\begin{equation}
\label{inter}
h_K^{s-1} \Vert \yyy - \yyy_I
\Vert_{s,K} \le C \Vert \yyy \Vert_{1,K} \le C' \, , \quad s=0,1
\end{equation}
for all $K \in \mc{C}_h$. Using the Cauchy--Schwartz
inequality, the bound (\ref{inter}) with $s=0$ and the definition
(\ref{shearnorm}) it follows
\begin{eqnarray}
\label{zz-split}
(\q -\q_h, \yyy) & = & (\q -\q_h, \yyy-\yyy_I) + (\q -\q_h, \yyy_I)
\nonumber \\
& \le & C \Vert \q - \q_h \Vert_{-1,h} + (\q -\q_h, \yyy_I) \, .
\end{eqnarray}
Note that $\yyy_I$ is both in $\VVh$ and $\Vstar$\,; moreover $\LL
\yyy_I = \nol$ on each  element $K$ of $\mc{C}_{h}$. As a
consequence, using (\ref{metheq}), (\ref{disc-shear}), (\ref{shear})
and Theorem \ref{consistency}, it follows
\begin{align}
\label{zz-2} 
(\q -\q_h, \yyy_I)
\nonumber
& =  a(\bbb-\bbb_h,\yyy_I)
     + \langle [\nabla  w_h - \bbb_h)] \cdot\s,
     M_{ns}(\yyy_I)\rangle_{\FREE}
\nonumber \\
& =: T_1 + T_2 \, .
\end{align}
Due to the continuity of the bilinear form and using bound
(\ref{inter}) with $s=1$ it immediately follows
\begin{equation}
\label{zz-3}
 T_1   \le   C \Vert \bbb-\bbb_h \Vert_{1}  \le   C
\trinorm (w-w_h,\bbb-\bbb_h) \trinorm_h \, .
\end{equation}
Using first the Cauchy--Schwartz inequality, then the Agmon
inequality, finally the bound (\ref{inter}) with $s=1$, Lemma
\ref{inverse} and the definition (\ref{3bar-norm}), we get
\begin{eqnarray}
\label{zz-4}
 T_2 & \le & \Big( \sum_{E \in \FREEh} h_E^{-1} \Vert
\nabla w_h - \bbb_h) \Vert_{0,E}^2 \Big)^{1/2} \Big( \sum_{E \in
\FREEh} h_E \Vert M_{ns}(\yyy_I) \Vert_{0,E}^2 \Big)^{1/2}
\nonumber \\
& \le & \Big( \sum_{K \in \mc{C}_{h}} h_K^{-2} \Vert \nabla  w_h - \bbb_h) \Vert_{0,K}^2 \Big)^{1/2} \Vert \yyy_I \Vert_1 \nonumber \\
& \le & C \trinorm (w-w_h,\bbb-\bbb_h) \trinorm_h \, ,
\end{eqnarray}
where in the last inequality we implicitly used the relation $\nabla
w - \bbb=\nol$. Combining (\ref{zz-1}), (\ref{zz-split}) with
(\ref{zz-2}), (\ref{zz-3}) and (\ref{zz-4}),  it follows that
\begin{equation}
\label{zz-7}
\Vert \q -\q_h \Vert_{-1,*}
\le C \big( \Vert \q -\q_h \Vert_{-1,h} 
+ \trinorm (w-w_h,\bbb -\bbb_h) \trinorm_h \big) \, .
\end{equation}
Joining (\ref{zz-7}), (\ref{minusoneh}), and using Theorem
\ref{apriori} the proposition immediately follows.
\end{proof}

The regularity of the solution to the Kirchhoff plate problems for
 convex polygonal domains, with all the three main types of boundary conditions,
is very case dependent. We refer for example to the work
\cite{MR80} in which a rather complete study is accomplished.
Note that, if $f \in H^{-1}(\Omega)$, in most cases of interest, the
regularity condition $w \in H^3(\Omega)$ is indeed achieved.

Note further, that with classical duality arguments and technical calculations it
is possible to derive the error bound
\begin{eqnarray}
\Vert w - w_h \Vert_1 \le C h^{s+1} \Vert w \Vert_{s+2} \, ,
\end{eqnarray}
if the regularity estimate
\begin{eqnarray}
\Vert w \Vert_3 \le C \Vert f \Vert_{-1}
\end{eqnarray}
holds. Moreover, if $k \ge 2$ and the regularity estimate
\begin{eqnarray}
\Vert w \Vert_4 \le C \Vert f \Vert_{0}
\end{eqnarray}
is satisfied, then it holds
\begin{eqnarray}
\Vert w - w_h \Vert_0 \le C h^{s+2} \Vert w \Vert_{s+2} \, .
\end{eqnarray}

\section{A-posteriori error estimates}
\label{4}
\setcounter{equation}{0}

In this section, we prove the reliability and the efficiency for an
a-posteriori error estimator for our method. To this end, we
introduce
\begin{eqnarray}
\label{eta-i}
&& \tilde\eta_K^2 
:= h_K^4 \Vert f + \ddiv \q_h \Vert_{0,K}^2 
+ h_K^{-2} \Vert \nabla w_h - \bbb_h \Vert_{0,K}^2 \, , 
\\
\label{eta-e} 
&& \eta_\edge^2 
:= \he^3 \Vert \llbracket \q_h \cdot \nn \rrbracket \Vert_{0,\edge}^2 
+ \he \Vert \llbracket \MM (\bbb_h) \nn \rrbracket \Vert_{0,\edge}^2 \, , 
\\
\label{eta-s} 
&& \eta_{S,E}^2 
:= \he \Vert \mnn (\bbb_h) \Vert_{0,\edge}^2 \, , 
\\
\label{eta-f}
&& \eta_{F,E}^2 
:= \he \Vert \mnn (\bbb_h) \Vert_{0,\edge}^2 
+ \he^3 \Vert \frac{\partial}{\partial s} \mns(\bbb_h)
-\q_h\cdot\nn \Vert_{0,\edge}^2 \, , 
\end{eqnarray}
where  $\he$ denotes the length of the edge $\edge$ and $\llbracket
\cdot \rrbracket$ represents the jump operator (which is assumed to
be equal to the function value on boundary edges). Further, for a
triangle $K\in \mc{C}_h$ we denote the sets of edges lying in the interior 
of $\Omega$, on $\SIMP$ and on $\FREE$,
by $I(K)$, $S(K)$ and $F(K)$, respectively. By $\simpedges$
we denote the set of all edges on $\SIMP$ and by $ \intedges$ the
ones lying in the interior of the domain.

Given any element $K \in \mc{C}_h$, let the local error indicator be
\begin{eqnarray}
\label{eq:eta-K}
\eta_K := \Big( \tilde\eta_K^2 
+ \frac{1}{2} \! \sum_{\edge\in I(K) } \eta_\edge^2 
+ \sum_{\edge \in S(K)} \eta_{S,E}^2 
+ \sum_{\edge\in F(K)} \eta_{F,E}^2 \Big)^{1/2} \, ,
\end{eqnarray}
Finally, the global error indicator is defined as
\begin{eqnarray}
\label{global} 
\eta := \Big( \sum_{K \in \mc{C}_h} \eta_K^2
\Big)^{1/2} \, .
\end{eqnarray}

\begin{remark}
It is worth noting that, by the definition (\ref{disc-shear}),
\begin{eqnarray}
\label{qplusL}
(\q_h + \LL \bbb_h)_{|K} 
= \frac{1}{\alpha h_K^2} (\nabla w_h - \bbb_h)_{|K} 
\qquad \forall K \in \mc{C}_h \, ,
\end{eqnarray}
which is the reason why there appears no terms of the kind $\Vert
\q_h + \LL \bbb_h \Vert_{0,K}$ in the error estimator. 
We note as well that scaling arguments give
\begin{eqnarray}
\sum_{\edge \in \FREEh} \he^{-1} \Vert \nabla w_h - \bbb_h \Vert_{0,\edge}^2 
\le C \sum_{K \in \mc{C}_h} h_K^{-2} \Vert \nabla w_h - \bbb_h \Vert_{0,K}^2 \, ,
\end{eqnarray}
which is the reason why there appears no boundary terms of the kind
$ \Vert \nabla w_h - \bbb_h \Vert_{0,\edge} $.
\end{remark}

\subsection{Upper bound}

In order to derive the reliability of the method we need the following saturation assumption.

\begin{assu}
Given a mesh $\mc{C}_h$, let $\mc{C}_{h/2}$ be the mesh obtained by
splitting each triangle $K \in \mc{C}_h$ into four triangles
connecting the edge midpoints. Let $(w_{h/2},\bbb_{h/2})$ be the
discrete solution corresponding to the mesh $\mc{C}_{h/2}$. We
assume that there exists a constant $ \rho $, $ 0<\rho<1 $, such
that \begin{eqnarray} 
\label{sat2}
 \trinorm 
 && (w-w_{h/2},\bbb-\bbb_{h/2}) \trinorm_{h/2} 
 + \Vert \q - \q_{h/2} \Vert_{-1,*}
\nonumber \\
 &&
 \quad \le \rho \big( \trinorm (w-w_h,\bbb-\bbb_h) \trinorm_h 
 + \Vert \q -\q_h \Vert_{-1,*}  \big) \, ,
\end{eqnarray}
where by $\trinorm \cdot \trinorm_{h/2}$ we indicate the mesh
dependent norm with respect to the new mesh $\mc{C}_{h/2}$.
\end{assu}

In the sequel,  we will  need the following result.

\begin{lemma}
\label{interlemma}
Let, for $v \in W_{h/2}$, the local seminorm be
\begin{eqnarray}
\label{starnorm}
\vert v \vert_{2,h/2,K} = \big( \sum_{K'\in \mc{C}_{h/2}, \ K' \subset K}
\vert v \vert_{2,K'}^2 \big)^{1/2} \, .
\end{eqnarray}
Then, there is a positive constant $C$ such that for all $v \in
W_{h/2}$ there exists $v_I \in W_h$ with the bound
\begin{eqnarray}
\Vert v - v_I \Vert_{0,K} 
+ h_K^{1/2}\Vert v - v_I \Vert_{0,\partial K} 
\le C h_K^2 \vert v \vert_{2,h/2,K} 
\qquad \forall K \in \mc{C}_h \, .
\end{eqnarray}
Moreover, $v_I$ interpolates $v$ at all the vertices of the triangulation $\mc{C}_{h/2}$.
\end{lemma}

\begin{proof}
We choose $v_I$ as the only function in $H^1(\Omega)$ such that
\begin{eqnarray}
&& v_{I |K} \in P_2(K) \qquad  \forall K \in \mc{C}_h \, ,
\nonumber \\
&& v_I(x) = v(x) \qquad  \forall x \in \mc{V}_{h/2} \, ,
\end{eqnarray}
where $\mc{V}_{h/2}$ represents the set of all the vertices of $\mc{C}_{h/2}$.
Note that it is trivial to check that $v_I \in W_h$ for all $k\ge 1$.
Observing that
\begin{eqnarray}
\vert v \vert_{2,h/2,K} 
+ \sum_{x\in \mc{V}_{h/2} \cap K} \vert v(x) \vert \, , 
\qquad v \in W_{h/2} \, , K \in \mc{C}_h \, ,
\end{eqnarray}
is indeed a norm on the finite dimensional space of 
the functions $v \in W_{h/2}$ restricted to $K$, 
the result follows applying the classical scaling argument.
\end{proof}

For simplicity, in the sequel we will treat the case $\SIMP =
\emptyset$, the general case following with identical arguments as
the ones that follow. We have the following preliminary result.

\begin{theorem}
\label{preliminary}
It holds
\begin{eqnarray}
\trinorm (w_{h/2}-w_h,\bbb_{h/2}-\bbb_h) \trinorm_{h/2}  
\le C \eta \, .
\end{eqnarray}
\end{theorem}

\begin{proof}
{\em Step 1.} Due to the stability of the discrete formulation,
proved in Theorem \ref{stability}, there exists a couple $(v,\yyy)
\in W_{h/2} \times \VV\!_{h/2}$ such that
\begin{eqnarray}
\label{stabound}
\trinorm (v,\yyy) \trinorm_{h/2} \le C
\end{eqnarray}
and
\begin{eqnarray}
\label{XX}
\trinorm (w_{h/2} - w_h , \bbb_{h/2} - \bbb_h) \trinorm_{h/2} \le
\mc{A}_{h/2} (w_{h/2} - w_h, \bbb_{h/2} - \bbb_h ; v, \yyy) \, .
\end{eqnarray}
Furthermore, we  have
\begin{eqnarray}
\label{X1}
\mc{A}_{h/2} (w_{h/2}, \bbb_{h/2} ; v, \yyy)  = (f,v) \, .
\end{eqnarray}

{\em Step 2.}  Simple calculations and the definition
(\ref{disc-shear}) give
\begin{eqnarray}
\label{X2}
&& \mc{B}_{h/2} (w_{h}, \bbb_{h} ; v, \yyy)
= a(\bbb_h, \yyy) - \sum_{K \in \mc{ C }_{h/2}} \alpha h_K^2 (\LL \bbb_h, \LL \yyy)_K
\nonumber  \\
&& \quad + \sum_{K \in \mc{ C }_{h/2}} \frac{1}{\alpha h_K^2} (\nabla w_h - \bbb_h - \alpha h_K^2 \LL \bbb_h,
    \nabla v - \yyy - \alpha h_K^2 \LL \yyy)_K
\nonumber \\
&& = a(\bbb_h, \yyy)
     - \sum_{K \in \mc{ C }_{h/2}} (\nabla w_h - \bbb_h , \LL \yyy)_K
     + \sum_{K \in \mc{ C }_{h/2}} (\q_h , \nabla v - \yyy)_K
\nonumber \\
&& \quad + R_1(w_h,\bbb_h; v,\yyy)
\nonumber \\
&& = \mc{B}_{h} (w_{h}, \bbb_{h} ; v, \yyy) + R_1(w_h,\bbb_h; v,\yyy) \, ,
\end{eqnarray}
where $ \q_h $ is defined as in (\ref{disc-shear}), i.e.,
based on the coarser mesh, and
\begin{align}
\label{Rdef}
R_1(w_h,\bbb_h; v,\yyy) 
= & \sum_{K \in \mc{ C }_{h/2}} \frac{1}{\alpha h_K^2} 
(\nabla w_h - \bbb_h , \nabla v - \yyy)_K \nonumber 
\\
& - \sum_{K \in \mc{ C }_h} \frac{1}{\alpha h_K^2} 
(\nabla w_h - \bbb_h , \nabla v - \yyy)_K \, .
\end{align}
The last term on the right hand side is well defined since $\nabla v
-\yyy$ is piecewise $L^2$-regular.

Let now  $\halfFREEh$   indicate   the set of all edges of
$\mc{C}_{h/2}$ lying   on $\FREE$. Adding and subtracting the
difference between the two forms it then follows
\begin{eqnarray}
\label{X4} 
\mc{D}_{h/2}(w_h,\bbb_h; v, \yyy) =
\mc{D}_{h}(w_h,\bbb_h; v, \yyy) + R_2(w_h, \bbb_h; v, \yyy) \, ,
\end{eqnarray}
where
\begin{align}
\label{Rpdef} 
R_2(w_h, \bbb_h; v, \yyy) = &
\sum_{\edge \in \halfFREEh} \frac{\gamma}{\he} 
\langle[\nabla w_h - \bbb_h ]\cdot\s , [\nabla v - \yyy]\cdot \s\rangle_\edge 
\nonumber \\
& - \sum_{\edge \in \FREEh} \frac{\gamma}{\he} 
\langle[\nabla w_h - \bbb_h ]\cdot\s , [\nabla v - \yyy]\cdot \s\rangle_\edge \, ,
\end{align}
and where the first member on the right hand side is indeed well
defined due to the piecewise regularity of $(v,\yyy)$. We will
denote
\begin{equation}
R(w_h, \bbb_h; v, \yyy)
= R_1(w_h, \bbb_h; v, \yyy) + R_2(w_h, \bbb_h; v, \yyy) \, .
\end{equation}
Joining (\ref{X1})--(\ref{Rpdef}) then yields
\begin{equation}
\label{rs1}
\mc{A}_{h /2} (w_{h }, \bbb_{h } ; v, \yyy)
= \mc{A}_{h} (w_{h }, \bbb_{h } ; v, \yyy)+  R (w_h, \bbb_h; v, \yyy) \, .
\end{equation}

{\em Step 3}. Let $v_I\in W_h$ be the interpolant defined in Lemma
\ref{interlemma} and let $\yyy_I\in \VVh$ be the piecewise linear
interpolant to $\yyy $. First, we have
\begin{equation}
\mc{A}_{h} (w_{h}, \bbb_{h} ; v_I, \yyy_I) = (f,v_I) \, .
\end{equation}
This, together with  (\ref{X1}) and (\ref{rs1}) gives
\begin{eqnarray}
\nonumber
\label{ww}
& &
\mc{A}_{h/2} (w_{h/2} - w_h, \bbb_{h/2} - \bbb_h ; v, \yyy)
\\
& &\quad = \mc{A}_{h/2} (w_{h/2} , \bbb_{h/2}  ; v, \yyy) 
    - \mc{A}_{h/2} ( w_h, \bbb_h ; v, \yyy)
\nonumber \\
&& \quad  = \mc{A}_{h/2} (w_{h/2} , \bbb_{h/2}  ; v, \yyy)
    - \mc{A}_{h} (w_{h }, \bbb_{h } ; v, \yyy)
    - R (w_h, \bbb_h; v, \yyy) 
\nonumber  \\
&& \quad =(f ,v-v_I ) 
    - \mc{A}_{h} (w_{h }, \bbb_{h } ;  v-v_I, \yyy-\yyy_I)
    - R (w_h, \bbb_h; v, \yyy) \, .
\end{eqnarray}

{\em Step 4.} Next, we bound  the last terms above.
Recalling that $\mc{ C }_{h/2}$ is a
subdivision of $\mc{ C }_{h}$, the Cauchy--Schwartz inequality,
(\ref{3bar-norm}) and  (\ref{stabound}) give
\begin{eqnarray}
\label{Rbound} 
&& \vert R_1(w_h,\bbb_h; v,\yyy) \vert 
\le 2 \vert \sum_{K \in \mc{ C }_{h/2}}
 \frac{1}{\alpha h_K^2} (\nabla w_h - \bbb_h , \nabla v - \yyy)_K \vert
\nonumber \\
&& \quad \le 2 \Big( \sum_{K \in \mc{ C }_{h/2}} \frac{1}{h_K^2} \Vert
\nabla w_h - \bbb_h \Vert_{0,K}^2 \Big)^{1/2} \Big( \sum_{K \in \mc{
C }_{h/2}} \frac{1}{h_K^2} \Vert \nabla v - \yyy \Vert_{0,K}^2 \Big)^{1/2}
\nonumber \\
&& \quad \le C \Big( \sum_{K \in \mc{ C }_{h/2}} \frac{1}{h_K^2} \Vert
\nabla w_h - \bbb_h \Vert_{0,K}^2 \Big)^{1/2} \, .
\end{eqnarray}
Using scaling  and arguments similar to those already adopted in
(\ref{Rbound}) it can be checked that
\begin{eqnarray}
\label{Rpbound} 
\vert R_2(w_h, \bbb_h; v, \yyy) \vert 
\le C \big( \sum_{K \in \mc{ C }_{h/2}} \frac{1}{h_K^2} 
\Vert \nabla w_h - \bbb_h \Vert_{0,K}^2 \big)^{1/2} \, .
\end{eqnarray}
Combining (\ref{Rbound}) and (\ref{Rpbound}) we get
\begin{eqnarray}
\label{RRbound}
\vert R(w_h,\bbb_h; v,\yyy) \vert 
\le \vert R_1(w_h,\bbb_h; v,\yyy) \vert 
+ \vert R_2(w_h,\bbb_h; v,\yyy) \vert
\le C \eta \, .
\end{eqnarray}

{\em Step 5.} Next, we expand, substitute the expression
(\ref{disc-shear}) for $\q_h$ and regroup the terms:
\begin{eqnarray}
\label{ff}
\nonumber
\lefteqn{(f ,v-v_I ) -
\mc{A}_{h} (w_{h }, \bbb_{h } ;  v-v_I, \yyy-\yyy_I)}
\\&
=& (f ,v-v_I ) -\Big\{a(\bbb_h, \yyy-\yyy_I) - \sum_{K \in \mc{ C
}_{h}} \alpha h_K^2 \big(\LL \bbb_h, \LL (\yyy-\yyy_I)\big)_K
\nonumber  \\
&& \quad + \sum_{K \in \mc{ C }_{h}} \frac{1}{\alpha h_K^2}
\big(\nabla w_h - \bbb_h - \alpha h_K^2 \LL \bbb_h,
    \nabla (v-v_I) - (\yyy-\yyy_I) - \alpha h_K^2 \LL( \yyy-\yyy_I)\big)_K
\nonumber \\
& & \quad + \big\langle m_{ns}(\bbb_h),[\nabla (v-v_I) -
(\yyy-\yyy_I)]\cdot \s\big\rangle_{\FREE} + \big\langle[\nabla w_h -
\bbb_h]\cdot\s, m_{ns}(\yyy-\yyy_I)\big\rangle_{\FREE}
\nonumber \\
&& \quad +  \sum_{\edge \in \FREEh}\frac{\gamma}{h_\edge}
\big\langle [\nabla w_h - \bbb_h]\cdot\s , [\nabla (v-v_I) -
(\yyy-\yyy_I)]\cdot \s \big\rangle_\edge\Big\} \nonumber
\\
& =& (f ,v-v_I )-\Big\{a(\bbb_h, \yyy-\yyy_I) - \sum_{K \in \mc{ C
}_{h}} \alpha h_K^2 \big(\LL \bbb_h+\q_h, \LL (\yyy-\yyy_I)\big)_K
\nonumber  \\
&& \quad +  \big(\q_h,
    \nabla (v-v_I) - (\yyy-\yyy_I)  \big)
\nonumber \\
& & \quad + \big\langle m_{ns}(\bbb_h),[\nabla (v-v_I) -
(\yyy-\yyy_I)]\cdot \s\big\rangle_{\FREE} + \big\langle[\nabla w_h -
\bbb_h]\cdot\s, m_{ns}(\yyy-\yyy_I)\big\rangle_{\FREE}
\nonumber \\
&& \quad +  \sum_{\edge \in \FREEh}\frac{\gamma}{h_\edge}
\big\langle [\nabla w_h - \bbb_h]\cdot\s , [\nabla (v-v_I) -
(\yyy-\yyy_I)]\cdot \s \big\rangle_\edge\Big\} \nonumber
\\
& =& \Big\{(f ,v-v_I )-
   \big(\q_h, \nabla (v-v_I) \big) 
   - \big\langle m_{ns}(\bbb_h),[\nabla (v-v_I) ]\cdot \s\big\rangle_{\FREE}
\nonumber \\
&& \qquad -  \sum_{\edge \in \FREEh}\frac{\gamma}{h_\edge}
\big\langle [\nabla w_h - \bbb_h]\cdot\s , [\nabla (v-v_I)  ]\cdot
\s \big\rangle_\edge\Big\} \nonumber
\\
&  &\quad - \Big\{a(\bbb_h, \yyy-\yyy_I) - \sum_{K \in \mc{ C }_{h}}
\alpha h_K^2 \big(\LL \bbb_h+\q_h, \LL (\yyy-\yyy_I)\big)_K-
\big(\q_h, \yyy-\yyy_I \big)
\nonumber \\
& & \qquad - \big\langle m_{ns}(\bbb_h),[\yyy-\yyy_I]\cdot
\s\big\rangle_{\FREE} + \big\langle[\nabla w_h - \bbb_h]\cdot\s,
m_{ns}(\yyy-\yyy_I)\big\rangle_{\FREE}
\nonumber \\
&& \qquad -  \sum_{\edge \in \FREEh}\frac{\gamma}{h_\edge}
\big\langle [\nabla w_h - \bbb_h]\cdot\s , [\yyy-\yyy_I]\cdot
\s \big\rangle_\edge\Big\} \nonumber
\nonumber \\
&=:& A - B \, . 
\end{eqnarray}

{\em Step 6.} In the part $A$ above, integration by parts, and
using the fact that $v(x)=v_I(x)$ at the corner points $x\in
\vertexset$, yields
\begin{eqnarray}
\label{A1}
& & (f ,v-v_I )- \big(\q_h, \nabla (v-v_I) \big) 
   - \big\langle m_{ns}(\bbb_h),[\nabla (v-v_I) ]\cdot \s\big\rangle_{\FREE}
\nonumber \\
& & \quad = (f+\ddiv \q_h ,v-v_I )
   + \big\langle \frac{\partial}{\partial \s}m_{ns}(\bbb_h)-\q_h\cdot\nn,
     v-v_I \big\rangle_{\FREE} \, .
\end{eqnarray}
The separate terms are then estimated as follows, using the Cauchy--Schwartz
inequality and Lemma \ref{interlemma},
\begin{eqnarray}
 \big\vert  \lefteqn{ (f+\ddiv \q_h ,v-v_I )  \big\vert
 =  \big\vert\sum_{K \in \mc{ C }_h} (f_h + \ddiv
\q_h, v - v_I)_K  \big\vert }
\nonumber \\
&& \quad
\le \Big( \sum_{K \in \mc{ C }_h} h_K^{4} \Vert f  + \ddiv
\q_h \Vert_{0,K}^2 \Big)^{1/2} \Big( \sum_{K \in \mc{ C }_h}
h_K^{-4} \Vert v-v_I \Vert_{0,K}^2 \Big)^{1/2}
\nonumber \\
&& \quad
\le C \Big( \sum_{K \in \mc{ C }_h} h_K^{4} \Vert f + \ddiv
\q_h \Vert_{0,K} \Big)^{1/2} \Big( \sum_{K \in \mc{ C }_{h}} \vert v
\vert_{2,h/2,K}^2 \Big)^{1/2}
\nonumber \\
&& \quad\le C \Big( \sum_{K \in \mc{ C }_h} \tilde\eta_K^2
\Big)^{1/2}
\end{eqnarray}
and
\begin{eqnarray}
\lefteqn{  \big\vert \big\langle \frac{\partial}{\partial
\s}m_{ns}(\bbb_h)-\q_h\cdot\nn, v-v_I
 \big\rangle_{\FREE}  \big\vert
=  \big\vert\sum_{\edge \in \FREEh} \langle
\frac{\partial}{\partial \s} \mns(\bbb_h) - \q_h \cdot \nn , v -
v_I\rangle_\edge    \big \vert } 
\nonumber \\
&& \quad \le \Big( \sum_{\edge \in \FREEh} \he^3 \Vert
\frac{\partial}{\partial \s} \mns(\bbb_h) - \q_h \cdot \nn
\Vert_{0,\edge}^2 \Big)^{1/2} \Big( \sum_{\edge \in \FREEh} \he^{-3}
\Vert v-v_I \Vert_{0,\edge}^2 \Big)^{1/2}
\nonumber \\
&& \quad \le C \Big( \sum_{\edge \in \FREEh} \eta_{F,E}^2
\Big)^{1/2} 
       \Big( \sum_{K \in \mc{C}_h} \vert v \vert_{2,h/2,K}^2 \Big)^{1/2}
\nonumber \\
&& \quad \le C \Big( \sum_{\edge \in \FREEh} \eta_{F,E}^2 \Big)^{1/2} \, .
\end{eqnarray}
The last term in $A$ is readily estimated by scaling estimates and
Lemma \ref{interlemma}
\label{A2}
\begin{eqnarray}
\lefteqn{\big\vert \sum_{\edge \in \FREEh}\frac{\gamma}{h_\edge}
\big\langle [\nabla w_h - \bbb_h]\cdot\s , [\nabla (v-v_I)  ]\cdot
\s \big\rangle_\edge\big\vert }
\nonumber \\
 && \quad \leq
\Big( \sum_{\edge \in \FREEh} \he^{-1}  \Vert
\nabla w_h - \bbb_h \Vert_{0,\edge}^2 \Big)^{1/2} \Big( \sum_{\edge
\in \FREEh} \he^{-1}  \Vert \nabla (v-v_I)  \Vert_{0,\edge}^2
\Big)^{1/2}
\nonumber \\
 && \quad \leq C \Big( \sum_{K \in \mc{C}_h} h_K^{-2}  \Vert
\nabla w_h - \bbb_h \Vert_{0,K}^2 \Big)^{1/2} ( \sum_{\edge \in
\FREEh} \he^{-3}  \Vert   v-v_I  \Vert_{0,\edge}^2 \Big)^{1/2}
\nonumber
\\
&& \quad  \leq C \Big( \sum_{K \in \mc{ C }_h} \tilde\eta_K^2 \Big)^{1/2} 
\Big( \sum_{K \in \mc{C}_h} \vert v \vert_{2,h/2,K}^2 \Big)^{1/2} 
\leq C \Big( \sum_{K \in \mc{ C }_h} \tilde\eta_K^2 \Big)^{1/2} \, . 
\end{eqnarray}
Collecting (\ref{A1})--(\ref{A2}) we obtain
\begin{equation}
\label{aest}
\big\vert A\big\vert \leq C \eta \, .
\end{equation}

{\em Step 7}. We will now estimate the term $B$. The following terms
are directly estimated as the similar terms above
\begin{eqnarray}
\label{B1}
&& \Big \vert  \big\langle[\nabla w_h - \bbb_h]\cdot\s,
m_{ns}(\yyy-\yyy_I)\big\rangle_{\FREE}\Big \vert + \Big\vert
\sum_{\edge \in \FREEh}\frac{\gamma}{h_\edge} \big\langle [\nabla
w_h - \bbb_h]\cdot\s , [\yyy-\yyy_I] \cdot \s
\big\rangle_\edge\Big\vert
\nonumber \\
&& \quad \leq C \eta \, .
\end{eqnarray}
Since $\yyy_I$ is piecewise linear, it holds $\LL \yyy_{I \vert K}=\nol$.
The inverse estimate then gives
\begin{eqnarray}
& & \Big\vert \sum_{K \in \mc{ C }_{h}} \alpha h_K^2
\big(\LL \bbb_h+\q_h, \LL (\yyy-\yyy_I)\big)_K\Big\vert 
= \Big\vert \sum_{K \in \mc{ C }_{h}} \alpha h_K^2 
\big(\LL \bbb_h+\q_h, \LL  \yyy \big)_K\Big\vert
\nonumber \\
& & \quad  \leq C \Big( \sum_{K \in \mc{ C }_{h}} \alpha h_K^2 
\Vert\LL \bbb_h+\q_h\Vert_{0,K}^2 \Big)^{1/2}\Vert \yyy \Vert_1
\nonumber \\
& & \quad \leq C \eta \, ,
\end{eqnarray}
where we in the last step used (\ref{qplusL}). The final step in
estimating the term $ B $ is to integrate by parts, use the Cauchy--Schwartz
inequality, interpolation estimates and again (\ref{qplusL}):
\begin{eqnarray}
\label{B2}
\nonumber 
& & \big\vert a(\bbb_h, \yyy-\yyy_I) 
- \big(\q_h, \yyy-\yyy_I \big) 
- \big\langle m_{ns}(\bbb_h),[\yyy-\yyy_I]\cdot \s\big\rangle_{\FREE}\big\vert
\nonumber \\
& & \quad = \Big\vert -\sum_{K \in \mc{C}_h}(\LL\bbb_h+\q_h,
\yyy-\yyy_I) + \sum_{E\in\intedges}\langle\llbracket \MM (\bbb_h)
\nn \rrbracket , \yyy-\yyy_I\rangle_E
\nonumber \\
& & \qquad \quad + \langle m_{nn}(\bbb_h),[\yyy-\yyy_I]\cdot\nn\rangle_{\SIMP\cup\FREE}
\Big\vert
\nonumber
\\
& & \quad \leq   \sum_{K \in \mc{C}_h}\Vert\LL\bbb_h+\q_h
\Vert_{0,K} \Vert  \yyy-\yyy_I\Vert_{0,K} +
\sum_{E\in\intedges}\Vert \llbracket \MM (\bbb_h) \nn \rrbracket
\Vert_{0,E}\Vert  \yyy-\yyy_I \Vert_{0,E} 
\nonumber
\\
& & \qquad  \quad + \sum_{E\in\simpedges\cup \FREEh}  \Vert
m_{nn}(\bbb_h) \Vert_{0,E} \Vert  \yyy-\yyy_I \Vert_{0,E}
\nonumber \\
&& \quad \leq C \eta \, . 
\end{eqnarray}
Collecting (\ref{B1})--(\ref{B2}) we obtain
\begin{equation}
\label{best}
\big\vert B\big\vert \leq C \eta \, .
\end{equation}

{\em Step 8.} The asserted estimate now follows from (\ref{XX}),
(\ref{ww}), (\ref{RRbound}), (\ref{ff}), (\ref{aest}) and (\ref{best}).
\end{proof}

We also have the following lemma for the shear force:

\begin{lemma}
\label{shearlemma}
It holds 
\begin{equation} 
\Vert \q_{h/2} -\q_h \Vert_{-1,*} \le C
\big( \trinorm (w_{h/2}-w_h,\bbb_{h/2}-\bbb_h) \trinorm_{h/2} + \eta \big) \, .
\end{equation}
\end{lemma}

\begin{proof}
We start by observing that, referring to the definition
(\ref{disc-shear}) and its "$h/2$" counterpart, $\q_h$ and
$\q_{h/2}$ are defined on different meshes and therefore with
different $h_K^2$ coefficients. However, recalling that the size
ratio between the two meshes is bounded, it is easy to check that an
opportune splitting and the triangle inequality give
\begin{eqnarray}
\label{zz-5} 
&& \Vert \q_{h/2} -\q_h \Vert_{-1,h}^2 
\le C \Big( \sum_{K \in \mc{C}_{h/2}} 
\Vert \nabla (w_{h/2} - w_h) - (\bbb_{h/2} - \bbb_h) \Vert_{0,K}^2
\nonumber \\
&& \quad + \sum_{K \in \mc{C}_{h}} \Vert \nabla w_h - \bbb_h \Vert_{0,K}^2
+ \sum_{K \in \mc{C}_{h/2}} h_K^2 \Vert \LL\bbb_{h/2} - \LL\bbb_h
\Vert_{0,K}^2 \Big) \, .
\end{eqnarray}
The first and the last term in
(\ref{zz-5}) can be bounded in terms of the $ \trinorm \cdot
\trinorm_{h/2} $ norm, simply using the definition (\ref{3bar-norm})
and the inverse inequality
\begin{equation}
h_K^2 \Vert \LL\bbb_{h/2} - \LL\bbb_h \Vert_{0,K}^2 \le C \Vert
\bbb_{h/2} - \bbb_h \Vert_{1,K}^2 \, .
\end{equation}
Therefore, recalling the definition (\ref{eta-i}), we get
\begin{eqnarray}
\label{zz-6}
&& \Vert \q_{h/2} -\q_h \Vert_{-1,h} 
\le C \big(\trinorm (w_{h/2}-w_h,\bbb_{h/2}-\bbb_h) \trinorm_{h/2} + \eta\big) \, .
\end{eqnarray}
The transition from the $\Vert \q_{h/2} -\q_h \Vert_{-1,h}$ norm to
the $\Vert \q_{h/2} -\q_h \Vert_{-1,*}$ norm is accomplished by using the
"Pitk\"{a}ranta--Verf\"{u}rth trick" with steps almost identical to
those used in Lemma \ref{shearlemma2}, therefore omitted.
\end{proof}

Joining Theorem \ref{preliminary} and Lemma \ref{shearlemma} gives
the following a-posteriori upper bound for the method.

\begin{theorem}
It holds 
\begin{equation} 
\trinorm (w-w_h,\bbb-\bbb_h) \trinorm_h + \Vert \q -\q_{h} \Vert_{-1,*}
\: \le \: C  \eta \, .
\end{equation}
\end{theorem}

\begin{proof}
Theorem \ref{preliminary} combined with Lemma \ref{shearlemma}
trivially gives
\begin{equation}
\label{sr}
\trinorm (w_{h/2}-w_h,\bbb_{h/2}-\bbb_h) \trinorm_{h/2} 
+ \Vert \q_{h/2} - \q_h \Vert_{-1,*}
\le C \eta \, .
\end{equation}
From the saturation assumption it follows
\begin{eqnarray}
& & \trinorm (w -w_h,\bbb  -\bbb_h) \trinorm_{h/2} 
+ \Vert \q - \q_h \Vert_{-1,*} 
\nonumber \\
& & \quad \leq \frac{1}{1-\rho} 
\big( \trinorm (w_{h/2}-w_h,\bbb_{h/2}-\bbb_h) \trinorm_{h/2} 
+ \Vert \q_{h/2} - \q_h \Vert_{-1,*}\big) 
\end{eqnarray}
and hence the assertion follows from (\ref{sr}).
\end{proof}

\subsection{Lower bound}

In this section, we prove the efficiency of the error estimator.
Given any edge $ \edge $ of the triangulation, we define $
\omega_\edge $ as the set of all the triangles $ K \in \mc{C}_h $
that have $ \edge $ as an edge. Given any $ K \in \mc{C}_h $, we
define $ \omega_K $ as the set of all the triangles in  $ \mc{C}_h $
that share an edge with $ K $. We then have the following lemma
\cite{LS+}:

\begin{lemma}
Given any edge $ \edge $ of the triangulation $ \mc{C}_h $, let $
P_k(\edge) $ be the space of polynomials of degree at most $ k $ on
$ \edge $. There exists a linear operator
\begin{equation}
    \Pi_\edge \: : \: P_k(\edge) \longrightarrow H^2_0(\omega_\edge)
\end{equation}
such that for all $ p_k \in P_k(\edge) $ it holds
\begin{align}
\label{Pi1}
    C_1 \Vert p_k \Vert_{0,\edge}^2
    & \le \langle p_k, \Pi_\edge(p_k) \rangle_\edge
    \le \Vert p_k \Vert_{0,\edge}^2 \, ,
    \\
\label{Pi2}
    \Vert \Pi_\edge (p_k) \Vert_{0,\omega_\edge}
    & \le C_2 \he^{1/2} \Vert p_k \Vert_{0,\edge} \, ,
\end{align}
where the positive constants $ C_i $ above depend only on $ k $ and
the minimum angle of the triangles in $ \mc{C}_h $.
\end{lemma}

Next, we define a local counterpart of the negative norm defined in
\eqref{eq:neg*norm} for the shear force.
\begin{equation}
\label{eq:neg*norm-local}
    \| \rr \|_{-1,*,\omega_K}
    = \sup_{\substack{ \yyy \in V_* \\ \yyy = \nol \ {\rm in} \ \Omega \backslash \omega_K}}
       \frac{ \la \rr, \yyy \ra }{ \| \yyy \|_1 } \, .
\end{equation}
We then have the following reliability result:

\begin{theorem}
It holds
\begin{equation}
\label{new}
    \eta_K
    \: \le \: C \big( \trinorm (w - w_h, \bbb - \bbb_h) \trinorm_{h,\omega_K}
    + \| \q - \q_h \|_{-1,*,\omega_K}
    + h_K^2 \| f - f_h \|_{0,\omega_K} \big) \, ,
\end{equation}
where $f_h$ is some approximation of the load $f$. Here $ \trinorm
\cdot \trinorm_{h,\omega_K} $ and $ \| \cdot \|_{0,\omega_K}$
represent, respectively, the standard restrictions of the norms $
\trinorm \cdot \trinorm_h $ and $ \| \cdot \|_0 $ to the domain $
\omega_K $.
\end{theorem}

\begin{proof}
The proof of the theorem consists of bounding separately all the
addenda of $ \eta_K $ in \eqref{eq:eta-K}.

{\em Step 1.}
We first bound the terms of $ \tilde\eta_K^2 $ in \eqref{eta-i}.
Considering the right hand side of \eqref{new}, the triangle
inequality immediately shows that it is sufficient to bound the term
$h_K^2 \Vert f_h + \ddiv \q_h \Vert_{0,K}$.

Given any $ K \in \mc{C}_h $, let $ b_K $ indicate the standard
third order polynomial bubble function on $ K $, scaled such that $
\| b_K \|_{L^\infty(K)} = 1 $. Given $ K \in \mc{C}_h $, let now $
\varphi_K \in H^2_0(K) $ be defined as
\begin{equation}
    \varphi_K = (f_h + \ddiv \q_h) \: b_K^2 \, .
\end{equation}
The standard scaling arguments then easily show that
\begin{align}
    \| f_h + \ddiv \q_h \|_{0,K}^2
    & \le C (f_h + \ddiv \q_h, \varphi_K)_K \, ,
    \\
\label{bubble2}
    \| \varphi_K \|_{0,K}
    & \le C \Vert f_h + \ddiv \q_h \Vert_{0,K} \, .
\end{align}

For the first term in $ \tilde\eta_K^2 $, the equilibrium equation
\eqref{strong2} and integration by parts give
\begin{eqnarray}
\label{Y1}
    && h_K^2 \Vert f_h + \ddiv \q_h \Vert_{0,K}^2
    \le C h_K^2 (f_h + \ddiv \q_h,\varphi_K)_K
\nonumber \\
    && \quad = C h_K^2 \big( (f + \ddiv \q_h, \varphi_K)_K + (f_h - f, \varphi_K)_K \big)
\nonumber \\
    && \quad = C h_K^2 \big( (-\ddiv \q + \ddiv \q_h, \varphi_K)_K + (f_h - f, \varphi_K)_K \big)
\nonumber \\
    && \quad = C h_K^2 \big( (\q_h - \q, \nabla \varphi_K)_K + (f_h - f, \varphi_K)_K \big) \, .
\end{eqnarray}
We note, in particular, that $ \nabla \varphi_K \in \Vstar$ and $
\nabla \varphi_K = \nol $ in $ \Omega \backslash K $. Therefore, the
duality inequality and the Cauchy--Schwartz inequality, followed by
the inverse inequality and the bound \eqref{bubble2} lead to the
estimate
\begin{eqnarray}
\label{Y2}
    && C h_K^2 \big( (\q_h - \q, \nabla \varphi_K)_K
       + (f_h - f, \varphi_K) \big)
\nonumber \\
    && \quad \le C \Vert \q - \q_h \Vert_{-1,*,K} \, h_K^2 \Vert \nabla \varphi_K \Vert_{1,K}
                 + C h_K^2 \Vert f - f_h \Vert_{0,K} \Vert \varphi_K \Vert_{0,K}
\nonumber \\
    && \quad \le C \big( \Vert \q - \q_h \Vert_{-1,*,K}
                 + h_K^2 \Vert f - f_h \Vert_{0,K} \big) \Vert f_h + \ddiv \q_h \Vert_{0,K} \, .
\end{eqnarray}
Combining now \eqref{Y1} with \eqref{Y2} gives
\begin{equation}
\label{hello}
    h_K^2 \Vert f_h + \ddiv \q_h \Vert_{0,K}
    \le C \big( \Vert \q - \q_h \Vert_{-1,*,K}
        + h_K^2 \Vert f - f_h \Vert_{0,K} \big) \, .
\end{equation}

The second term of $ \tilde\eta_K^2 $ in \eqref{eta-i} can be
directly bounded by using the Kirchhoff condition \eqref{strong3}
with the definitions \eqref{eq:h-norm}--\eqref{3bar-norm},
\begin{align}
    h_K^{-1} \Vert \nabla w_h - \bbb_h \Vert_{0,K}
    & = h_K^{-1} \Vert \nabla (w - w_h) - (\bbb - \bbb_h) \Vert_{0,K}^2
\nonumber \\
    & \le \trinorm (w - w_h, \bbb - \bbb_h) \trinorm_{h,K} \, .
\end{align}

{\em Step 2.}
We next bound the terms of $ \eta_\edge^2 $ in \eqref{eta-e}. Given
now $ \edge \in I(K) $, an edge of the the element $ K $ lying in
the interior of $ \Omega $, let
\begin{equation}
\label{Y5}
    \boldsymbol{\varphi}_\edge
    = \Pi_\edge(\llbracket \MM (\bbb_h) \nn \rrbracket) \, ,
\end{equation}
where, with a little abuse of notation, the operator $ \Pi_\edge $
is intended as applied on each single component. Then, from
\eqref{Pi1} with integration by parts, it follows that
\begin{align}
\label{Y3}
    \he^{1/2} \Vert \llbracket \MM (\bbb_h) \nn \rrbracket \Vert_{0,\edge}^2
    & \le C \he^{1/2} \langle \llbracket \MM (\bbb_h) \nn \rrbracket,
            \boldsymbol{\varphi}_\edge \rangle_\edge
\nonumber \\
    & = C \he^{1/2} \big( (\LL\bbb_h, \boldsymbol{\varphi}_\edge)_{\omega_\edge}
        + (\MM(\bbb_h), \nabla \boldsymbol{\varphi}_\edge)_{\omega_\edge} \big) \, ,
\end{align}
where we recall that $ \omega_\edge $ was defined at the start of
this section. Integration by parts and the equation \eqref{strong1}
immediately lead to the identity
\begin{equation}
    (\MM(\bbb), \nabla \boldsymbol{\varphi}_\edge)_{\omega_\edge}
    = - (\LL\bbb, \boldsymbol{\varphi}_\edge)_{\omega_\edge}
    = (\q, \boldsymbol{\varphi}_\edge)_{\omega_\edge} \, ,
\end{equation}
which, applied to \eqref{Y3}, gives
\begin{eqnarray}
\label{Y4}
    && \he^{1/2} \Vert \llbracket \MM (\bbb_h) \nn \rrbracket \Vert_{0,\edge}^2
\nonumber \\
    && \quad \le C \he^{1/2} \big( (\LL\bbb_h + \q, \boldsymbol{\varphi}_\edge)_{\omega_\edge}
                 + (\MM(\bbb_h) - \MM(\bbb),\nabla \boldsymbol{\varphi}_\edge)_{\omega_\edge} \big)
\nonumber \\
    && \quad = C \he^{1/2} \big( (\LL\bbb_h + \q_h, \boldsymbol{\varphi}_\edge)_{\omega_\edge}
               + (\q - \q_h,\boldsymbol{\varphi}_\edge)_{\omega_\edge}
\nonumber \\
    && \qquad + (\MM(\bbb_h) - \MM(\bbb), \nabla \boldsymbol{\varphi}_\edge)_{\omega_\edge} \big) \, .
\end{eqnarray}

Next, we bound the three terms on the right hand side of \eqref{Y4}.
For the first term, the identity \eqref{qplusL}, the
Cauchy--Schwartz inequality, the definition \eqref{Y5} and the bound
\eqref{Pi2} give
\begin{eqnarray}
\label{Y4-1}
    && \he^{1/2} (\LL\bbb_h + \q_h, \boldsymbol{\varphi}_\edge)_{\omega_\edge}
\nonumber \\
    && \quad \le C \Big( \sum_{K \subset \: \omega_\edge}
                 h_K^{-2} \Vert \nabla w_h - \bbb_h \Vert^2_{0,K} \Big)^{1/2}
                 \Vert \llbracket \MM (\bbb_h) \nn \rrbracket \Vert_{0,\edge}
\nonumber \\
    && \quad \le C \: \trinorm (w - w_h, \bbb - \bbb_h) \trinorm_{h,\omega_\edge}
                 \Vert \llbracket \MM (\bbb_h) \nn \rrbracket \Vert_{0,\edge} \, .
\end{eqnarray}
For the second term on the right hand side of \eqref{Y4}, we note
that $ \boldsymbol{\varphi}_\edge \in \Vstar$ and $
\boldsymbol{\varphi}_\edge = \nol $ in $ \Omega \backslash
\omega_\edge $. Therefore, the duality inequality and the definition
\eqref{Y5} combined with the bound \eqref{Pi2} give
\begin{align}
\label{Y4-2}
    \he^{1/2}(\q - \q_h, \boldsymbol{\varphi}_\edge)_{\omega_\edge}
    & \le  \he^{1/2} \Vert \q - \q_h \Vert_{-1,*,\omega_\edge}
           \Vert \boldsymbol{\varphi}_\edge \Vert_{1,\omega_\edge}
\nonumber \\
    & \le  C \Vert \q - \q_h \Vert_{-1,*,\omega_\edge}
           \Vert \llbracket \MM (\bbb_h) \nn \rrbracket \Vert_{0,\edge} \, .
\end{align}
For the third term of \eqref{Y4}, the Cauchy--Schwartz inequality,
then the inverse inequality and finally \eqref{Y5} combined with the
bound \eqref{Pi2} lead to the estimate
\begin{eqnarray}
\label{Y4-3}
    && \he^{1/2} (\MM(\bbb_h) - \MM(\bbb), \nabla \boldsymbol{\varphi}_\edge)_{\omega_\edge}
    \le  C \Vert \bbb - \bbb_h \Vert_{1,\omega_\edge}
         \: h_K^{-1/2} \Vert \boldsymbol{\varphi}_\edge \Vert_{0,\omega_\edge}
\nonumber \\
    && \quad \le  C \Vert \bbb - \bbb_h \Vert_{1,\omega_\edge}
                  \Vert \llbracket \MM (\bbb_h) \nn \rrbracket \Vert_{0,\edge} \, .
\end{eqnarray}
Now, by combining \eqref{Y4-1}, \eqref{Y4-2} and \eqref{Y4-3} with
\eqref{Y4} it follows
\begin{equation}
    \he^{1/2} \Vert \llbracket \MM (\bbb_h) \nn \rrbracket \Vert_{0,\edge}
    \le C \big( \trinorm (w - w_h, \bbb - \bbb_h) \trinorm_{h,\omega_\edge}
        + \Vert \q - \q_h \Vert_{-1,*,\omega_\edge} \big) \, .
\end{equation}

The remaining term of $ \eta_\edge^2 $ is bounded with similar
arguments; with the notation
\begin{equation}
\label{Y5-2}
    \varphi_\edge
    = \Pi_\edge(\llbracket \q_h \cdot \nn \rrbracket) \, ,
\end{equation}
the identity
\begin{equation}
    - (\ddiv \q, \varphi_\edge)_{\omega_\edge}
    = (\q, \nabla \varphi_\edge)_{\omega_\edge}
\end{equation}
with \eqref{Y1} implies
\begin{align}
    \he^{1/2} \| \llbracket \q \cdot \nn \rrbracket \|_{0,\edge}^2
    & \le  C \he^{1/2} \langle \llbracket \q \cdot \nn \rrbracket, \varphi_\edge \rangle_\edge
\nonumber \\
    & \le  C \he^{1/2} \big( (f - f_h, \varphi_\edge)_{\omega_\edge}
           + (\q_h - \q, \nabla \varphi_\edge)_{\omega_\edge} \big) \, .
\end{align}
Finally, we note that $ \nabla \varphi_\edge \in \Vstar$ and $
\nabla \varphi_\edge = \nol $ in $ \Omega \backslash \omega_\edge $.
Therefore,
\begin{equation}
     \he^{3/2} \Vert \llbracket \q_h \cdot \nn \rrbracket \Vert_{0,\edge}
     \le C \big( \Vert \q - \q_h \Vert_{-1,*,\omega_\edge}
         + h_K^2 \Vert f - f_h \Vert_{0,\omega_\edge} \big) \, .
\end{equation}

{\em Step 3.}
Third, we bound the only term of $ \eta_{S,E}^2 $ in \eqref{eta-s}
which appears in $ \eta_{F,E}^2 $ as well. Given now a triangulation
edge $ \edge $ in $ S(K) \cup F(K) $, let
\begin{equation}
    \varphi_\edge = \Pi_\edge (m_{nn}(\bbb_h)) \, .
\end{equation}
Due to \eqref{Pi1} and \eqref{L-op}, integration by parts gives
(here $ \NABLA $ denotes the tensor valued gradient applied to a
vector valued function)
\begin{eqnarray}
\label{Y7}
    && \he^{1/2} \Vert m_{nn}(\bbb_h) \Vert_{0,\edge}^2
    \le \he^{1/2} \langle m_{nn}(\bbb_h - \bbb), \varphi_\edge\rangle_\edge
\nonumber \\
    && \quad = \he^{1/2} \langle \MM_{n}(\bbb_h - \bbb), \varphi_\edge \nn\rangle_\edge
\nonumber \\
    && \quad = \he^{1/2} \big( (\MM (\bbb_h - \bbb), \NABLA (\varphi_\edge \nn))_{\omega_\edge}
         + (\LL(\bbb_h - \bbb), \varphi_\edge \nn)_{\omega_\edge} \big) \, ,
\end{eqnarray}
where $ \nn $ is, as usual, the chosen normal unit vector to $ \edge
$. For the first term, using the Cauchy--Schwartz inequality, then
the inverse inequality and finally the bound \eqref{Pi2} we easily
get
\begin{align}
\label{Y8}
    \he^{1/2} (\MM (\bbb_h-\bbb), \NABLA(\varphi_\edge \nn))_{\omega_\edge}
    & \le  \he^{1/2} \Vert \bbb - \bbb_h \Vert_{1,\omega_\edge}
           \Vert \NABLA(\varphi_\edge \nn) \Vert_{0,\omega_\edge}
\nonumber \\
    & \le  C \Vert \bbb - \bbb_h \Vert_{1,\omega_\edge} \Vert m_{nn}(\bbb_h) \Vert_{0,\edge} \, .
\end{align}
For the second term in \eqref{Y7}, recalling \eqref{strong1} we have
\begin{eqnarray}
\label{pip}
    && \he^{1/2} (\LL(\bbb_h - \bbb), \varphi_\edge \nn)_{\omega_\edge}
\nonumber \\
    && \quad = \he^{1/2} (\LL\bbb_h + \q_h, \varphi_\edge \nn)_{\omega_\edge}
             + \he^{1/2} (\q - \q_h, \varphi_\edge \nn)_{\omega_\edge} \, .
\end{eqnarray}
Observing now that $ \varphi_\edge \nn \in \Vstar$ and $
\varphi_\edge \nn = \nol $ in $ \Omega \backslash \omega_\edge $,
the two terms on the right hand side of \eqref{pip} can be bounded
with the same arguments used above, respectively, in \eqref{Y4-1}
and \eqref{Y4-2}. Omitting the details, we therefore get
\begin{eqnarray}
\label{Y9}
    && \he^{1/2} (\LL(\bbb_h - \bbb), \varphi_\edge \nn)_{\omega_\edge}
        \le C \big( \trinorm (w - w_h, \bbb - \bbb_h) \trinorm_{h,\omega_\edge}
\nonumber \\
    && \quad + \Vert \q - \q_h \Vert_{-1,*,\omega_\edge} \big)
               \Vert m_{nn}(\bbb_h) \Vert_{0,\edge} \, .
\end{eqnarray}
From \eqref{Y7}, \eqref{Y8} and \eqref{Y9} we get
\begin{equation}
\label{Y12}
    \he^{1/2} \Vert m_{nn}(\bbb_h) \Vert_{0,\edge}
    \le C \big( \trinorm (w - w_h, \bbb - \bbb_h) \trinorm_{h,\omega_\edge}
        + \Vert \q - \q_h \Vert_{-1,*,\omega_\edge} \big) \, .
\end{equation}

{\em Step 4.}
Finally, we bound the last term of $ \eta_{F,E}^2 $ in
\eqref{eta-f}. Given now a triangulation edge $ \edge $ in $ F(K) $, let
\begin{equation}
\label{def2}
    \varphi_\edge
    = \Pi_\edge(\frac{\partial}{\partial s} \mns(\bbb_h) - \q_h\cdot\nn) \, .
\end{equation}
Using \eqref{Pi1} and recalling \eqref{bc-free}, we obtain
\begin{eqnarray}
\label{Ystart}
    && \he^{3/2} \Vert \frac{\partial}{\partial s} \mns(\bbb_h) - \q_h\cdot\nn \Vert_{0,\edge}^2
\nonumber \\
    && \quad \le \he^{3/2} \big( \langle\frac{\partial}{\partial s} \mns(\bbb_h - \bbb),
                 \varphi_\edge\rangle_{\edge}
                 + \langle[\q - \q_h] \cdot \nn, \varphi_\edge\rangle_{\edge} \big) \, .
\end{eqnarray}
For the first term, integration by parts on the edge and simple
algebra give
\begin{eqnarray}
\label{Y11}
    && \he^{3/2} \langle\frac{\partial}{\partial s} \mns(\bbb_h - \bbb),
       \varphi_\edge\rangle_{\edge}
    = \he^{3/2} \langle\mns(\bbb - \bbb_h), \nabla \varphi_\edge \cdot \s \rangle_{\edge}
\nonumber \\
    && \quad = \he^{3/2} \big( \langle\MM(\bbb - \bbb_h)\nn,
               \nabla \varphi_\edge\rangle_{\edge}
               - \langle m_{nn}(\bbb - \bbb_h), \nabla \varphi_\edge \cdot \nn \rangle_{\edge} \big) \, .
\end{eqnarray}
Using again integration by parts, the first term in \eqref{Y11} can
be written as
\begin{eqnarray}
\label{Y10}
    && \he^{3/2} \langle \MM(\bbb - \bbb_h)\nn, \nabla \varphi_\edge\rangle_{\edge}
\nonumber \\
    && \quad = \he^{3/2} \big( \LL(\bbb - \bbb_h), \nabla \varphi_\edge)_{\omega_\edge}
               + \langle \MM(\bbb-\bbb_h),\NABLA\nabla\varphi_\edge)_{\omega_\edge} \big) \, .
\end{eqnarray}
The second term in \eqref{Ystart}, again due to integration by parts
and recalling \eqref{strong2}, is instead equivalent to
\begin{eqnarray}
    && \he^{3/2}  \langle [\q - \q_h] \cdot \nn, \varphi_\edge\rangle_{\edge}
    = \he^{3/2} \big (\q - \q_h, \nabla \varphi_\edge)_{\omega_\edge}
\nonumber \\
    && \quad - (f_h + \ddiv \q_h, \varphi_\edge)_{\omega_\edge}
             - (f - f_h, \varphi_\edge)_{\omega_\edge} \big) \, .
\end{eqnarray}
For the first term, due to \eqref{strong1} and \eqref{disc-shear},
we now have
\begin{eqnarray}
\label{Y-40}
    && \he^{3/2} (\q - \q_h, \nabla \varphi_\edge)_{\omega_\edge}
\nonumber \\
    && \quad = \he^{3/2} \big( \LL (\bbb_h - \bbb), \nabla \varphi_\edge)_{\omega_\edge}
               - \frac{1}{\alpha h_{\omega_\edge}^2}
               (\nabla w_h - \bbb_h, \nabla \varphi_\edge)_{\omega_\edge} \big) \, ,
\end{eqnarray}
where $ h_{\omega_\edge} $ is the size of the triangle $
\omega_\edge $. Combining all the identities from \eqref{Ystart} to
\eqref{Y-40} it follows that
\begin{eqnarray}
\label{Y-41}
    && \he^{3/2} \Vert \frac{\partial}{\partial s} \mns(\bbb_h) - \q_h\cdot\nn \Vert_{0,\edge}^2
\nonumber \\
    && \quad \le \he^{3/2} \big( (\MM(\bbb - \bbb_h),
                 \NABLA \nabla\varphi_\edge)_{\omega_\edge}
                 - (m_{nn}(\bbb - \bbb_h), \nabla \varphi_\edge \cdot \nn \rangle_{\edge}
\nonumber \\
    && \qquad - \frac{1}{\alpha h_{\omega_\edge}^2} (\nabla w_h - \bbb_h,
             \nabla \varphi_\edge)_{\omega_\edge}
             - (f_h + \ddiv \q_h, \varphi_\edge)_{\omega_\edge}
\nonumber \\
    && \qquad - (f - f_h, \varphi_\edge)_{\omega_\edge} \big) \, .
\end{eqnarray}
For the second term on the right hand side of \eqref{Y-41},
recalling \eqref{bc-free}, using the Cauchy--Schwartz inequality and
the bound \eqref{Y12} we have
\begin{eqnarray}
\label{Y15}
    && \he^{3/2} \langle m_{nn}(\bbb - \bbb_h), \nabla \varphi_\edge \cdot \nn \rangle_\edge
    \le \he^{1/2} \Vert m_{nn}(\bbb_h) \Vert_{0,\edge}
        \: \he \Vert \nabla\varphi_\edge \Vert_{0,\edge}
\nonumber \\
    && \quad \le C \big( \trinorm (w - w_h,\bbb - \bbb_h)  \trinorm_{h,\omega_\edge}
                 + \Vert \q - \q_h \Vert_{-1,*,\omega_\edge} \big)
                 \he \Vert \nabla\varphi_\edge \Vert_{0,\edge} \, ,
\end{eqnarray}
which, using the inverse inequality and the bound \eqref{Pi2}, gives
\begin{eqnarray}
\label{Y16}
    && \he^{3/2} \langle m_{nn}(\bbb - \bbb_h), \nabla \varphi_\edge \cdot \nn \rangle_\edge
    \le C \big( \trinorm (w - w_h, \bbb - \bbb_h) \trinorm_{h,\omega_\edge}
\nonumber \\
    && \quad + \Vert \q - \q_h \Vert_{-1,*,\omega_\edge} \big)
             \Vert \frac{\partial}{\partial s} \mns(\bbb_h) - \q_h\cdot\nn \Vert_{0,\edge} \, .
\end{eqnarray}
The remaining terms on the right hand side of \eqref{Y-41} can all
be bounded using the Cauchy--Schwartz inequality, the inverse
inequality and the bounds \eqref{hello}, \eqref{Pi2} as already
shown for the similar previous cases. Without showing all the
details, we finally get
\begin{eqnarray}
    && \he^{3/2} \Vert \frac{\partial}{\partial s} \mns(\bbb_h)
                 - \q_h\cdot\nn \Vert_{0,\edge}^2
\nonumber \\
    && \quad \le C \big( \trinorm (w - w_h, \bbb - \bbb_h) \trinorm_{h,\omega_\edge}
                 + \Vert \q - \q_h \Vert_{-1,*,\omega_\edge}
\nonumber \\
    && \qquad + h_K^2 \Vert f - f_h \Vert_{0,K} \big)
              \Vert \frac{\partial}{\partial s} \mns(\bbb_h)
              - \q_h \cdot \nn \Vert_{0,\edge} \, ,
\end{eqnarray}
or, trivially,
\begin{eqnarray}
    && \he^{3/2} \Vert \frac{\partial}{\partial s} \mns(\bbb_h)
       - \q_h\cdot\nn \Vert_{0,\edge}
\nonumber \\
    && \quad \le C \big( \trinorm (w - w_h, \bbb - \bbb_h) \trinorm_{h,\omega_\edge}
                 + \Vert \q - \q_h \Vert_{-1,*,\omega_\edge}
                 + h_K^2 \Vert f - f_h \Vert_{0,K} \big) \, .
\end{eqnarray}
Recalling now the definitions for $ \eta_K $ in \eqref{eta-i} and
the local negative norm in \eqref{eq:neg*norm-local}, the
proposition is proved.
\end{proof}


\bibliography{Khhoff}
\bibliographystyle{siam}

\end{document}